\documentclass[12pt, a4paper]{article}
\usepackage{supertabular, setspace,hyperref}
\usepackage{amsmath,amssymb,amsthm,textcomp}
\usepackage{amsfonts,graphicx}
\usepackage[mathscr]{eucal}
\usepackage{color}
\usepackage{float}
\usepackage[caption = false]{subfig}
\usepackage{tikz}
\usetikzlibrary
{shapes,arrows,chains,matrix,positioning,scopes,shadows,calc,snakes}
\theoremstyle{definition}

\numberwithin{equation}{section}

\newcommand{\ncom}{\newcommand}

\ncom{\beq}{\begin{equation}}
\ncom{\eeq}{\end{equation}}
\ncom{\bea}{\begin{eqnarray*}}
\ncom{\eea}{\end{eqnarray*}}
\ncom{\beqa}{\begin{eqnarray}}
\ncom{\eeqa}{\end{eqnarray}}
\ncom{\nno}{\nonumber}
\ncom{\non}{\nonumber}
\ncom{\ds}{\displaystyle}
\ncom{\half}{\frac{1}{2}}
\ncom{\mbx}{\makebox{.25cm}}
\ncom{\hs}{\mbox{\hspace{.25cm}}}
\ncom{\rar}{\rightarrow}
\ncom{\Rar}{\Rightarrow}
\ncom{\noin}{\noindent}
\ncom{\bc}{\begin{center}}
\ncom{\ec}{\end{center}}
\ncom{\sz}{\scriptsize}
\ncom{\rf}{\ref}
\ncom{\s}{\sqrt{2}}
\ncom{\sgm}{\sigma}
\ncom{\Sgm}{\Sigma}
\ncom{\psgm}{\sigma^{\prime}}
\ncom{\dt}{\delta}
\ncom{\Dt}{\Delta}
\ncom{\lmd}{\lambda}
\ncom{\Lmd}{\Lambda}
\ncom{\Th}{\Theta}
\ncom{\e}{\eta}
\ncom{\eps}{\epsilon}
\ncom{\pcc}{\stackrel{P}{>}}
\ncom{\lp}{\stackrel{L_{p}}{>}}
\ncom{\dist}{{\rm\,dist}}
\ncom{\sspan}{{\rm\,span}}
\ncom{\re}{{\rm Re\,}}
\ncom{\im}{{\rm Im\,}}
\ncom{\sgn}{{\rm sgn\,}}
\ncom{\ba}{\begin{array}}
\ncom{\ea}{\end{array}}
\ncom{\hone}{\mbox{\hspace{1em}}}
\ncom{\htwo}{\mbox{\hspace{2em}}}
\ncom{\hthree}{\mbox{\hspace{3em}}}
\ncom{\hfour}{\mbox{\hspace{4em}}}
\ncom{\vone}{\vskip 2ex}
\ncom{\vtwo}{\vskip 4ex}
\ncom{\vonee}{\vskip 1.5ex}
\ncom{\vthree}{\vskip 6ex}
\ncom{\vfour}{\vspace*{8ex}}
\ncom{\norm}{\|\;\;\|}
\ncom{\integ}[4]{\int_{#1}^{#2}\,{#3}\,d{#4}}
\ncom{\vspan}[1]{{{\rm\,span}\{ #1 \}}}
\ncom{\dm}[1]{ {\displaystyle{#1} } }
\ncom{\ri}[1]{{#1} \index{#1}}

\newtheoremstyle
   {remarkstyle}
   {}
   {11pt}
   {}
   {}
   {\bfseries}
   {:}
   {     }
   {\thmname{#1} \thmnumber{#2} }

\theoremstyle{remarkstyle}

\def\eps{\varepsilon}

\begin{document}
\title{\bf On shock reflection-diffraction in a van der Waals gas}
\author{\bf Neelam Gupta and V. D. Sharma\\
{\it Department of Mathematics,
Indian Institute of Technology Bombay,}\\
{\it Powai, Mumbai-400076}}
\date{}
\maketitle
\begin{abstract}
\noindent The problem of a weak shock, reflected and diffracted by a wedge, is studied for the two-dimensional compressible Euler system. Some recent developments are overviewed and a perspective is presented within the context of a real gas, modeled by the van der Waals equation of state. The regular reflection configuration and the detachment criterion are studied in the light of real gas effects. Some basic features of the phenomenon and the nature of the self-similar flow pattern are explored using asymptotic expansions. The analysis presented here predicts several inviscid flow properties of the real gases undergoing shock reflection-diffraction phenomenon.
\end{abstract}
\vone
\section{Introduction}
Shock reflection problem, which has captured the interest of researchers over the years, is one of the most important problems for the mathematical theory of multidimensional conservation laws that is still largely incomplete. The experimental, computational and asymptotic analyses show that various patterns of reflected shocks may occur including the regular and Mach reflections (see, Courant and Friedrichs \cite{courant1976}, Glimm and Majda \cite{glimm}, Glass \cite{glass}, Zheng \cite{zheng,zheng1}, Ben-Dor \cite{dor}, Chen \cite{chen}, and Chang \& Hsiao \cite{chang}).
When a weak plane shock hits a wedge head on, two processes take place simultaneously. The incident shock wave is reflected by the wedge surface and at the same time the flow behind it is deflected by the wedge corner, producing a nearly circular diffracted wave expanding from the vertex; the circular wave emanating from near the vertex moves at sonic speed of the incident flow. Here, we consider the case in which only regular reflection is expected to occur; for small disturbance approximation in weak shock reflection, this corresponds to relatively large wedge angles. The reflected shock, which is weak enough, travels backward to join the diffracted wave smoothly. The flow in the diffracted wave region was calculated by Keller and Blank \cite{keller1951} using the linearized theory; this solution was modified later by Hunter and Keller \cite{kel} using the theory of weakly nonlinear geometrical acoustics \cite{keller1983, majda}. They showed that the diffracted wavefronts of linear acoustics are actually shocks. An approximate analytical solution to the shock-wedge diffraction problem was proposed by Harabetian \cite{eduard1987} using an alternate method that invokes perturbation expansions with multiple scales. Much effort has been devoted to the study of this problem through simplified models that capture various features of Euler system within the context of an ideal gas (see, Morawetz \cite{morawetz}, Zheng \cite{zheng2}, Rosales and Tabak \cite{tabak}, Brio and Hunter \cite{Hunter2000}, Canic et al. \cite{keyfitz}, and Hunter \& Tesdall \cite{tesdall}). It is well known that at high pressure or low temperature, the behavior of gases deviates from the ideal gas-law and follows van der Waals type gas that deals with the possible real gas effects (without phase-transition); examples cover a family of shock wave problems with complicated interface patterns and a hydrodynamic model of sonoluminescene (an acoustic-induced light emission phenomenon) (see,\cite{ cramer, kluwick, Wu1996, arora, manoj, manoj1}). In the present paper, we study the regular reflection configuration and the detachment criterion when the real gas effects are taken into account. In the limit of vanishing van der Waals excluded volume, we recover the result obtained by Chang and Chen \cite{chen1986}, which is a refinement of Von Neumann's criterion. One of the main objectives of the present paper is to study how the real gas effects influence the behavior of the local structure of the self-similar solutions of the compressible Euler equations near a singular point; the motivation stems from the work carried out in \cite{keller1951, kel, eduard1987, chen1986}. The real gas effects, presented here, are characterized by a van der Waals type equation of state. The analysis presented here predicts several inviscid properties of real gases undergoing shock reflection-diffraction phenomenon; a summary of the results is presented in the last section.  \\
The set-up of the reflection consists of a straight shock hitting a wedge at the origin at time $t=0$. The shock is assumed to be weak and moving, parallel to the y-axis, towards the wedge which is placed symmetrically about the flow direction, namely the x-axis, (Figure $1(i)$). The gas ahead of the shock is at rest. Since the problem is symmetric with respect to the x-axis, it suffices to consider the problem in the upper half plane outside the half wedge.\\
The basic equations of this study are the Euler equations
\begin{align}\label{equ1}
\begin{split}
\rho_{t} + ({\rho}{u})_x + ({\rho}{v})_y ={0},\\  
({\rho}{u})_{t} + ({\rho}{u^2}+p)_x + ({\rho}{u}{v})_y ={0},\\
({\rho}{v})_{t} + ({\rho}{u}{v})_x + ({\rho}{v^2}+p)_y ={0},\\
\left({\rho}{\Big(e+\frac{{u^2}+{v^2}}{2}\Big)}\right)_{t} + \left({\rho}{u}{\Big(h+\frac{{u^2}+{v^2}}{2}\Big)}\right)_{x}+\left({\rho}{v}{\Big(h+\frac{{u^2}+{v^2}}{2}\Big)}\right)_{y}=0,
\end{split}
\end{align}
where $\rho$, ($u$, $v$), $p$, $e$, and $h$ denote respectively the density, velocity components, pressure, internal energy, and specific enthalpy, while  $e$ and $h$ are given functions of $\rho$ and $p$ which satisfy the thermodynamical constraints $Td{\mbox{\scriptsize{S}}}=de+pdV=dh-Vdp$ with $T(\rho, p)$ being the temperature, $V$ the specific volume and ${\mbox{\scriptsize{S}}}(\rho, p)$ the specific entropy. We consider the situation when the gas obeys a van der Waals equations of state of the form   
\begin{equation}\label{equ}
p=\frac{RT}{(V-b)},~~~e=\frac{p(V-b)}{\gamma-1},~~~{\mbox{\scriptsize{S}}}=c_v\ln{(p(V-b)^{\gamma})}+constant,~~~h=\frac{p(\gamma{V}-b)}{\gamma-1},
\end{equation} 
where R is the gas constant, $\gamma(>1)$ the ratio of specific heats, and $b$ the van der Waals excluded volume. 
Consider a weak shock hitting a wedge with half angle $\alpha\in(0,\pi/2)$. The state ahead of the shock is ${(\rho, u, v, p)}={(\rho_0, 0, 0,  p_0)}$ for some $p_0>0$. The state behind the shock is ${( \rho_1, u_1, 0, p_1)}$ with $p_1>p_0$. So we seek a solution of the system (\ref{equ1}) with initial data
\begin{equation}\label{equ2}
{(\rho, u, v, p)}\Big|_{t=0} = \left\{
  \begin{array}{l l}
    {(\rho_0, 0, 0,  p_0)}, &\left|y\right|> x\tan{\alpha}, x>0\\
    {( \rho_1, u_1, 0, p_1)}, & x<0,\\
\end{array} \right.
\end{equation}
and the slip boundary condition along the wedge
\begin{equation}\label{equ3}
v=u{\tan\alpha}\Big|_{y=x{\tan\alpha}}\hspace{1cm} x>0,\hspace{.2cm} t>0.
\end{equation}
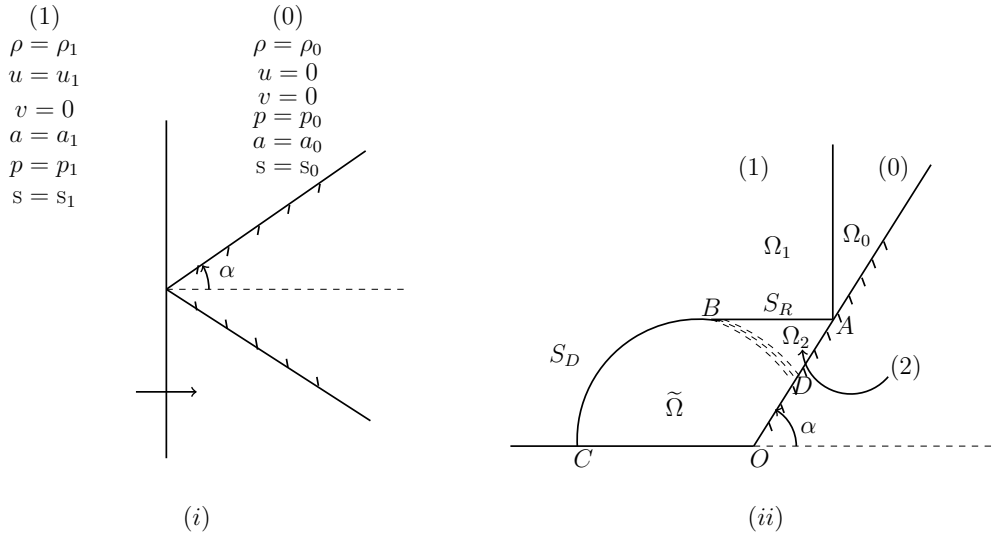
\begin{figure}[h]
\centering
\scalebox{.8}{
\begin{tikzpicture}
			\begin{scope}
			\draw [thick](0,.2)--(0,5.8);
			\draw [dashed](0,3)--(4,3);
			\draw [thick][rotate around={35:(0,3)}](0,3)--(4,3);
			\draw [thick][rotate around={-33:(0,3)}](0,3)--(4,3);
			\draw [->][thick](-0.5,1.3)--(0.5,1.3);
			\draw [->] [thick](0.7,3) arc (0:30:0.8); 
			\draw [thick][rotate around={10:(1,2.35)}](1,2.35)--(1,2.5);
			\draw [thick][rotate around={10:(1.5,2)}](1.5,2)--(1.5,2.2);
			\draw [thick][rotate around={10:(0.5,2.65)}](0.5,2.65)--(0.5,2.8);
			\draw [thick][rotate around={10:(2,1.7)}](2,1.7)--(2,1.85);
			\draw [thick][rotate around={10:(2.5,1.4)}](2.5,1.4)--(2.5,1.55);
			\draw [thick][rotate around={-10:(1,3.55)}](1,3.55)--(1,3.7);
			\draw [thick][rotate around={-10:(1.5,3.9)}](1.5,3.9)--(1.5,4.05);
			\draw [thick][rotate around={-10:(0.5,3.25)}](0.5,3.25)--(0.5,3.4);
			\draw [thick][rotate around={-10:(2,4.25)}](2,4.25)--(2,4.4);
			\draw [thick][rotate around={-10:(2.5,4.6)}](2.5,4.6)--(2.5,4.75);	
			\end{scope}
			\draw (1,3.3) node {$\alpha$};
			\draw (-2,7.5) node {$(1)$};
			\draw (-2,7) node {$\rho=\rho_{1}$};
			\draw (-2,6.5) node {$u=u_{1}$};
			\draw (-2,6) node {$v=0$};
			\draw (-2,5.5) node {$a=a_{1}$};
			\draw (-2,5.0) node {$p=p_{1}$};
			\draw (-2,4.5) node {$\mbox{\scriptsize{S}}=\mbox{\scriptsize{S}}_{1}$};
			\draw (2,7.5) node {$(0)$};
			\draw (2,7) node {$\rho=\rho_{0}$};
			\draw (2,6.6) node {$u=0$};
			\draw (2,6.2) node {$v=0$};
			\draw (2,5.8) node {$p = p_{0}$};
 			\draw (2,5.4) node {$a=a_{0}$};
			\draw (2,5.0) node {$\mbox{\scriptsize{S}}=\mbox{\scriptsize{S}}_{0}$};	
				\draw (.5,-.8) node {$(i)$};
\end{tikzpicture}\hspace{1.5cm}
\begin{tikzpicture}
			\begin{scope}
			\draw [thick](-4,3)--(0,3);
			\draw [dashed](0,3)--(4,3);
			\draw [thick][rotate around={58:(0,3)}](0,3)--(5.5,3);
			\draw [->] [thick](0.7,3) arc (0:60:0.7); 
			\draw [thick][rotate around={20:(0.3,3.25)}](0.3,3.25)--(0.3,3.4);
			\draw [thick][rotate around={20:(.5,3.55)}](.5,3.55)--(.5,3.7);
			\draw [thick][rotate around={20:(.7,3.85)}](.7,3.85)--(.7,4);
			\draw [thick][rotate around={20:(.87,4.15)}](.87,4.15)--(.87,4.3);
			
			\draw [thick][rotate around={20:(1.05,4.45)}](1.05,4.45)--(1.05,4.6);
			\draw [thick][rotate around={20:(1.26,4.75)}](1.26,4.75)--(1.26,4.9);
			\draw [thick][rotate around={20:(1.44,5.05)}](1.44,5.05)--(1.44,5.2);
			\draw [thick][rotate around={20:(1.62,5.35)}](1.62,5.35)--(1.62,5.5);
			\draw [thick][rotate around={20:(1.83,5.65)}](1.83,5.65)--(1.83,5.8);
			\draw [thick][rotate around={20:(2,5.95)}](2,5.95)--(2,6.1);
			\draw [thick][rotate around={20:(2.19,6.25)}](2.19,6.25)--(2.19,6.4);
			\draw [thick](1.3,5.1)--(1.3,8.0);
			\draw [thick](1.3,5.1)--(-0.7,5.1);
			\draw [-] [thick](-2.9,3) arc (3:-100:-2);
			\draw [dashed][rotate around={65:(-0.7,5.1)}](-0.7,5.1) arc (5:-30:2.75);
			\draw [dashed][rotate around={65:(-0.6,5.1)}](-0.6,5.1) arc (5:-30:2.7);
			\draw [dashed][rotate around={65:(-0.5,5.1)}](-0.5,5.1) arc (5:-30:2.7);
			\draw [<-] [thick](.8,4.6) arc (4.6:140:-0.8);
			\end{scope}
			\draw (0,7.6) node {$(1)$};
			\draw (0.4,6.3) node {$\Omega_1$};
			\draw (-1.3,3.7) node {$\widetilde{\Omega}$};
			\draw (2.3,7.6) node {$(0)$};
			\draw (1.7,6.5) node {$\Omega_0$};	
			\draw (0.7,4.8) node {$\Omega_2$};
			\draw (2.5,4.3) node {$(2)$};
			\draw (1.5,5.0) node {$A$};
			\draw (-0.7,5.3) node {$B$};
			\draw (0.4,5.35) node {$S_R$};
			\draw (-2.8,2.8) node {$C$};
			\draw (-3.1,4.5) node {$S_D$};
			\draw (.1,2.8) node {$O$};
			\draw (0.8,4) node {$D$};
			\draw (.9,3.3) node {$\alpha$};
			\draw (.2,1.8) node {$(ii)$};
\end{tikzpicture}
}
\caption{\textit{(i)~The primary shock hitting a wedge and (ii)~the self-similar form of the flow pattern restricted to the upper half-plane.}}
\label{figure1}
\end{figure}
\section{Self-similar flow and shock reflection-diffraction configuration}
Since the coefficients in equations governing the physical process do not depend on variables, equations (\ref{equ1}) together with initial and boundary conditions (\ref{equ2}) and (\ref{equ3}) are invariant under the dilation $t\rightarrow\nu{t}$, $x\rightarrow\nu{x}$, $y\rightarrow\nu{y}$, where $\nu>0$ is an arbitrary constant, and so we look for the solution with the property\\i.e.,
 \begin{equation*}
{(\rho, u, v, p)}(t, x, y)={(\rho, u, v, p)}(\nu{t}, \nu{x}, \nu{y}).
\end{equation*}
Thus, by introducing the variables ~$\zeta=\frac{\sqrt{{x^2}+{y^2}}}{t}$~ and $\theta=tan^{-1}\left(\frac{y}{x}\right)$, we express equations (\ref{equ1}) in self-similar polar form
\begin{align}\label{equ5}
\begin{split}
&\left({\rho(\mbox{\scriptsize{U}}-\zeta)}\right)_{\zeta}+\left({\frac{\rho\mbox{\scriptsize{V}}}{\zeta}}\right)_{\theta}+{\frac{\rho\mbox{\scriptsize{U}}}{\zeta}}+{\rho}=0,\\&
\left({\rho(\mbox{\scriptsize{U}}-\zeta)^2+p}\right)_{\zeta}+\left(\frac{\rho{(\mbox{\scriptsize{U}}-\zeta)}\mbox{\scriptsize{V}}}{\zeta}\right)_\theta+{\frac{\rho}{\zeta}}\left((\mbox{\scriptsize{U}}-\zeta)^2-\mbox{\scriptsize{V}}^2\right)+3{\rho}(\mbox{\scriptsize{U}}-\zeta)=0,\\&
\left({\rho}(\mbox{\scriptsize{U}}-\zeta)\mbox{\scriptsize{V}}\right)_\zeta+\left(\frac{{\rho}\mbox{\scriptsize{V}}^2+p}{\zeta}\right)_\theta+\frac{2\rho(\mbox{\scriptsize{U}}-\zeta)\mbox{\scriptsize{V}}}{\zeta}+3{\rho}\mbox{\scriptsize{V}}=0,\\&
\left(\rho(\mbox{\scriptsize{U}}-\zeta)\left(h+\frac{\mbox{\scriptsize{U}}^2+\mbox{\scriptsize{V}}^2}{2}\right)+\zeta{p}\right)_{\zeta}+\left(\frac{\rho\mbox{\scriptsize{V}}}{\zeta}\left(h+\frac{\mbox{\scriptsize{U}}^2+\mbox{\scriptsize{V}}^2}{2}\right)\right)_{\theta}+\rho{\left(h+\frac{\mbox{\scriptsize{U}}^2+\mbox{\scriptsize{V}}^2}{2}\right)}{\left(1+\frac{\mbox{\scriptsize{U}}}{\zeta}\right)}-p=0,
\end{split}
\end{align}
with initial and boundary conditions:
\begin{equation}\label{equ6}
 \lim_{{\zeta} \to \infty}{(\rho, u, v, p)}
  = \left\{
  \begin{array}{l l}
    {(\rho_0, 0, 0, p_0)}, & {\alpha}\leq{\theta}<{\frac{\pi}{2}},\\
		{(\rho_1, u_1, 0, p_1)}, & {\frac{\pi}{2}}\leq\theta<{\pi}.\\
\end{array} \right.
\end{equation}
and 
\begin{equation}\label{equ7}
v=u\tan{\theta}\Big|_{\theta=\alpha},
\end{equation}
where $h$ is given by $(\ref{equ})_4$, and  
\begin{equation}\label{equ5*}
\mbox{\scriptsize{U}}=u\cos{\theta}+v\sin{\theta},~~ \mbox{\scriptsize{V}}=-u\sin{\theta}+v\cos{\theta}. 
\end{equation}
If $\rho$ and ${\mbox{\scriptsize{S}}}$ are chosen as independent  variables, many calculations for Euler system (\ref{equ1}) can be simplified; for instance the speed of sound is 
\begin{equation}
a(\rho, {\mbox{\scriptsize{S}}})=\sqrt{\partial{p}/\partial{\rho}}=\sqrt{\frac{{\gamma}p}{\rho(1-b\rho)}};~~0\leq{b\rho}<1,
\end{equation}
 and the energy equation $(\ref{equ1})_4$, for smooth solutions, may be written as 
\begin{equation}\label{equS}
(\mbox{\scriptsize{U}}-\zeta)\mbox{\scriptsize{S}}_{\zeta}+({{\mbox{\scriptsize{V}}}/{\zeta}})\mbox{\scriptsize{S}}_{\theta}=0.
\end{equation}
It may be noticed that the unsteady Euler system (\ref{equ1}), governing the gas flow, is hyperbolic; however, the corresponding pseudo-stationary flow in self-similar coordinates is governed by mixed type equations. Indeed, the system (\ref{equ5}) changes its type from elliptic to hyperbolic when the point $(\zeta, \theta)$ runs from the origin to infinity. Since the problem is symmetric with respect to the line $\theta=0$, it suffices to consider the problem in the upper half plane, $\Omega=\left\{(\zeta,\theta):\zeta>0,\alpha\leq\theta\leq\pi/2\right\}$, outside the half wedge. As the shock front hits the wedge head on and propagates further along the wedge, it is reflected by the wedge surface at $A$ (see Figure $1(ii)$), whereas the induced flow behind the incident shock wave is diffracted by the wedge corner at sonic speed. Assuming that the shock reflection is a regular one, the location of the incident shock after it has moved beyond the domain of the influence of the origin (vertex of the wedge) is given by
\begin{equation}\label{equ7**}
\zeta=a_0\sec{\theta},
\end{equation}
where $a_0=\sqrt{{{\gamma}p_0}/{\rho_0(1-b\rho_0)}}$;
 the straight line segment $AB$ of the reflected shock is given by 
\begin{equation}\label{equ7*}
\zeta={a_0}\tan{\alpha}/(\sin(\theta-\alpha)\sec{\alpha}+\sin(2\alpha-\theta))\equiv{\zeta}^{*}.
\end{equation}
The unknown curved portion $BC$, due to the influence of the origin, joins the diffracted wavefront smoothly and gives rise the overall shock reflection-diffraction phenomenon. Here $BD$ is the fixed boundary referred to as the sonic arc, $\zeta=a_0$, across which there is a continuous transition from the supersonic region $\Omega_2$ to the subsonic region $\tilde{\Omega}$, whereas $BC$ is the free boundary, called diffraction of the planar shock, across which the transition undergoes a jump from the supersonic region ${\Omega}_1$ to the subsonic region $\tilde{\Omega}$ near the origin. These wavefronts, referred to as boundaries, separate the upper half $(\zeta, \theta)$-plane into four regions.
\begin{align}\nonumber
\begin{split}
&\Omega_0=\left\{(\zeta, \theta): \zeta>{a_0}\sec{\theta},~~ \alpha<\theta<{\pi}/2\right\},\\&
\Omega_1=\left\{(\zeta, \theta): \zeta^{*}<\zeta<{a_0}\sec{\theta},~~ \alpha<\theta<2\alpha\right\}\cup \left\{(\zeta, \theta): \zeta>{a_0},~~ 2\alpha<\theta<{\pi}\right\},\\&
\Omega_2=\left\{(\zeta, \theta): \zeta^{*}>\zeta>{a_0},~~ \alpha<\theta<2\alpha\right\},~~\widetilde{\Omega}=\left\{(\zeta, \theta): \zeta<{a_0},~~ \alpha<\theta<{\pi}\right\}.
\end{split}
\end{align}
From equation (\ref{equ7*}), it may be noticed that ${d{\zeta}^*}/{d\tilde{b}}>0$, where $\tilde{b}=b\rho_0$; this implies that an increase in  $\tilde{b}$ causes the domain $\tilde{\Omega}$ of linearized solution to become larger as compared to the ideal gas case.\\
Hence in order to determine the entire flow field and the wave structure, one needs to solve the free boundary value problem for a degenerate elliptic equation.
The system (\ref{equ5}) has four eigenvalues 
\begin{equation}\label{equ9}
\lambda=\frac{\mbox{\scriptsize{V}}}{\zeta(\mbox{\scriptsize{U}}-\zeta)},~~~\text{(multiplicity-2)}
\end{equation}
\begin{equation}\label{equ10}
\lambda=\frac{\mbox{\scriptsize{V}}(\mbox{\scriptsize{U}}-\zeta){\pm}a{\sqrt{\mbox{\scriptsize{V}}^2-a^2+(\mbox{\scriptsize{U}}-\zeta)^2}}}{{\zeta}\left((\mbox{\scriptsize{U}}-\zeta)^2-a^2\right)},
\end{equation}
with $\mbox{\scriptsize{V}}^2+(\mbox{\scriptsize{U}}-\zeta)^2>a^2$. Equations (\ref{equ10}) show that the system (\ref{equ5}) is hyperbolic with four eigenvalues and the flow is supersonic; when $\mbox{\scriptsize{V}}^2+(\mbox{\scriptsize{U}}-\zeta)^2<a^2$, the system is mixed type as two equations in (\ref{equ5}) are hyperbolic and the other two are elliptic. However, $\mbox{\scriptsize{V}}^2+(\mbox{\scriptsize{U}}-\zeta)^2=a^2$ represents a sonic curve in $(\zeta, \theta)$ plane. In general, the system (\ref{equ5}) is mixed type and the flow is transonic.
\section{State behind the incident and reflected shocks}
In order to find the states behind the incident and reflected shocks, denoted by subscripts-1 and -2 respectively, we need the Rankine-Hugoniot (R-H) conditions in 2D. Let $\zeta=\zeta(\theta)$ be a shock curve with slope ${\zeta}'{(\theta)}$. Then it follows from (\ref{equ5}) that
\begin{align}\label{equ11}
\begin{split}
&\left[{\rho}{(\mbox{\scriptsize{U}}-\zeta)}\right]{\zeta}d{\theta}=\left[\rho{\mbox{\scriptsize{V}}}\right]d{\zeta},\\&
[\rho(\mbox{\scriptsize{U}}-\zeta)^2+p]{\zeta}d{\theta}=[{\rho}(\mbox{\scriptsize{U}}-\zeta)\mbox{\scriptsize{V}}]d{\zeta},\\&
[{\rho}(\mbox{\scriptsize{U}}-\zeta)\mbox{\scriptsize{V}}]{\zeta}d{\theta}=[\rho{\mbox{\scriptsize{V}}^2}+p]d{\zeta},\\&
[\rho(\mbox{\scriptsize{U}}-\zeta)(h+\frac{\mbox{\scriptsize{U}}^2+\mbox{\scriptsize{V}}^2}{2})+{\zeta}p]{\zeta}d{\theta}=[\rho{\mbox{\scriptsize{V}}}(h+\frac{\mbox{\scriptsize{U}}^2+\mbox{\scriptsize{V}}^2}{2})]d{\zeta},
\end{split}
\end{align} 
where square brackets, [.], denote jumps across the shock. Let $q_t=(\mbox{\scriptsize{U}}-\zeta)d{\zeta}+\mbox{\scriptsize{V}}{\zeta}d{\theta}$, ~and~ $q_n=(\mbox{\scriptsize{U}}-\zeta){\zeta}d{\theta}-{\mbox{\scriptsize{V}}}d{\zeta}$ (with ${\zeta^2}d{\theta^2}+d{\zeta}^2=1$) be the components of pseudo-velocity vector $(\mbox{\scriptsize{U}}-\zeta, \mbox{\scriptsize{V}})$ along the tangent and normal to the shock curve $\zeta=\zeta(\theta)$. Then equations (\ref{equ11}) may be written as
\begin{equation}\label{equ12}
\left[{\rho}{q_n}\right]=0,~[{\rho}{q_t}{q_n}]=0,~[\rho{q_n}^2+p]=0,~\left[\rho{q_n}\left(h+\frac{{q_t}^2+{q_n}^2}{2}\right)\right]=0.
\end{equation}
When $q_{n0}\neq0$, it follows from $(\ref{equ12})_1$ that $q_{n1}\neq0$, showing thereby that the pseudo-flow is no longer tangential to the shock curve; indeed, it moves from state-0 to state-1 satisfying the entropy condition ${\rho_0}<{\rho_1}$. Let $M_0={q_0}/{a_0}$ (respectively, $M_1={q_1}/{a_1}$) be the up-stream (respectively, down-stream) shock Mach number relative to the up-stream (respectively, down-stream) flow.
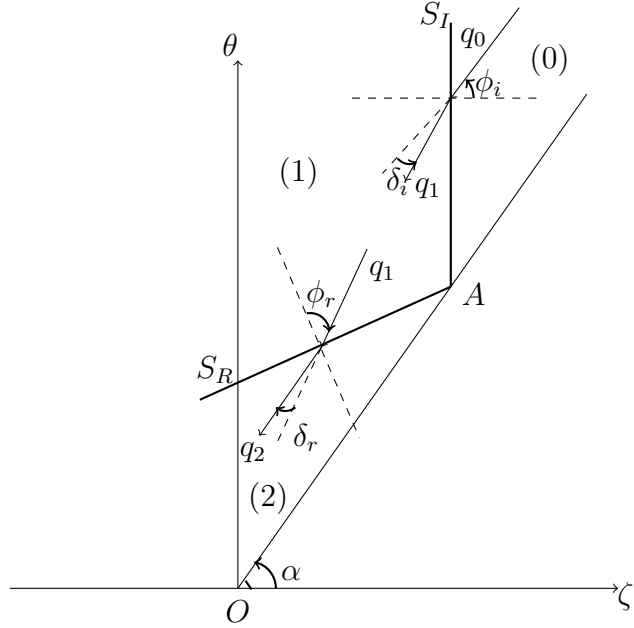
\begin{figure}[h]
\centering
\begin{tikzpicture}
			\begin{scope}
			\draw [->](-4,-4)--(4,-4);
	  \draw [->](-1,-4)--(-1,3);
			\draw [][rotate around={55:(-1,-4)}](-1,-4)--(7,-4);
   		\draw [thick](1.8,3.5)--(1.8,0);
  		\draw [thick](1.8,0)--(-1.5,-1.5);
  		\draw [dashed](.5,2.5)--(3,2.5);
  		\draw [->](2.7,3.7)--(1.8,2.5);
  		\draw [dashed](1.8,2.5)--(.9,1.5);
  		\draw [->](1.8,2.5)--(1.2,1.4);
  		\draw [dashed][rotate around={80:(.1,-.8)}](1.3,0)--(-1,-1.5);
 		  \draw [thick][rotate around={260:(-.9,-3.9)}](-.9,-3.9)--(-.8,-3.8);
 		  \draw [->](.7,.5)--(.1,-.8);
 		  \draw [dashed](.1,-.8)--(-.5,-2.1);
 		  \draw [->][rotate around={-10:(.1,-.8)}](.1,-.8)--(-.5,-2.1);
 		  \draw [->] [thick][](2.1,2.5) arc (0:45:.35);
	    \draw [->] [thick][rotate around={230:(1.1,1.7)}](1.1,1.7) arc (-10:70:.2);	
	    \draw [<-] [thick][](.2,-.6) arc (0:100:.25);
	    \draw [<-] [thick][rotate around={230:(-.5,-1.6)}](-.5,-1.6) arc (0:80:.2);
	    \draw [->] [thick][](-.5,-4) arc (0:80:.35);
		  \end{scope}
	\draw (4.1,-4.1) node {$\zeta$};
	\draw (-1.1,3.2) node {$\theta$};
		\draw (2.1,-.1) node {$A$};
		\draw (-1,-4.3) node {$O$};
	\draw (-1.3,-1.1) node{$S_R$};
	\draw (1.6,3.6) node {$S_{I}$};
		\draw (3.1,3) node {$(0)$};
	\draw (-.2,1.5) node {$(1)$};
	\draw (-.6,-2.8) node {$(2)$};
	\draw (2.1,3.3) node {$q_0$};
	\draw (1.5,1.3) node {$q_1$};
	\draw (.9,.2) node {$q_1$};
	\draw (-.8,-2.2) node {$q_2$};
	\draw (2.3,2.7) node {$\phi_i$};
	\draw (1.1,1.4) node {$\delta_i$};
	\draw (.1,-.1) node {$\phi_r$};	
	\draw (-.1,-2) node {$\delta_r$};	
	\draw (-.3,-3.8) node {$\alpha$};		
  \end{tikzpicture}
  \caption{\textit{The wave configuration of a regular reflection in a pseudo-steady flow.}}
  \label{figure1}
  \end{figure} 
Then, following \cite{chang}, it can be shown that in a van der Waals gas, the states $(\rho_0, 0, 0, p_0)$ and $(\rho_1, u_1, 0, p_1)$ on the two sides of the incident shock are related as 
\begin{equation}\label{equ21}
\dfrac{p_{1}}{p_0}=\dfrac{(\gamma+1-2\tilde{b})\beta_{i}-(\gamma-1)}{(\gamma+1)-(\gamma-1+2\tilde{b})\beta_{i}},~~~~\tan{\delta_i}=\dfrac{({\beta_i}-1)\tan{\phi_i}}{1+{\beta_i}{\tan}^2{\phi_i}},
\end{equation}
\begin{equation}\label{equM}
M^2_0=\dfrac{q^2_{n0}(1+{\tan}^2{\phi_i})}{a^2_0}=\dfrac{2{\beta_i}(1-\tilde{b})\sec^2{\phi_i}}{(\gamma+1)-{\beta_i}(\gamma-1+2\tilde{b})},~~~~~~~~~~~~
\end{equation}
\begin{equation}\label{equ23}
M^2_1=\dfrac{q^2_{n1}}{a_0^2}{\sec^2(\phi_i+\delta_i)}=\dfrac{2(1-\tilde{b}\beta_i)(1+{\beta^2_i}{\tan^2{\phi_i}})}{(\gamma+1){\beta_i}-(\gamma-1+2\tilde{b}\beta_i)},~~~~~~~~~~
\end{equation}
where $\beta_i={\rho_1}/{\rho_0}$, $\phi_i$ is the angle between the shock velocity vector $\stackrel{\rightarrow}{q_0}=(q_{t0}, q_{n0})$ and the shock normal with $0<\phi_i<\pi/2$, and $\delta_i$ is the angle between $\stackrel{\rightarrow}{q_0}$ and $\stackrel{\rightarrow}{q_1}$.\\
Further, since $p_1>p_0>0$ and  $0\leq\tilde{b}<1$, equation (\ref{equM}) implies that 
\begin{equation}\label{equ26*}
1<{\beta_i}<(\gamma+1)/(\gamma-1+2\tilde{b}).
\end{equation}
Similarly, the R-H conditions for the reflected shock can be written as 
\begin{equation}\label{equ24}
\dfrac{p_2}{p_1}=\dfrac{(\gamma+1-2\tilde{b}\beta_i){\beta_r}-(\gamma-1)}{(\gamma+1)-{\beta_r}(\gamma-1+2\tilde{b}\beta_i)},~~~M^2_1=\dfrac{2{\beta_r}(1+\tan^2{\phi_r})(1-\tilde{b}\beta_i)}{(\gamma+1)-{\beta_r}(\gamma-1+2\tilde{b}\beta_i)},
\end{equation}
\begin{equation}\label{equ26}
M^2_2=\dfrac{2(1+{\beta^2_r}\tan^2{\phi_r})(1-\tilde{b}{\beta_i}\beta_r)}{(\gamma+1){\beta_r}-(\gamma-1+2\tilde{b}{\beta_i}\beta_r)},~~~\tan{\delta_r}=\dfrac{({\beta_r}-1)\tan{\phi_r}}{1+{\beta_r}\tan^2{\phi}_r},
\end{equation}
\begin{equation}\label{equ28}
1<{\beta_r}<(\gamma+1)/(\gamma-1+2\tilde{b}\beta_i),
\end{equation}
where $\beta_r={\rho_2}/{\rho_1}$, $\phi_r$ is the angle between the shock velocity vector $\stackrel{\rightarrow}{q_1}=(q_{t1}, q_{n1})$ and the shock normal with $0<\phi_r<\pi/2$, and $\delta_r$ is the angle between $\stackrel{\rightarrow}{q_1}$ and $\stackrel{\rightarrow}{q_2}$.\\
The boundary condition (\ref{equ7}) requires that the state $(\rho_2, u_2, v_2, p_2)$ be such that $v_2=u_2{\tan{\alpha}}$; this condition together with R-H conditions determine the state-2. Let $\sigma_i=\tan{\delta_i}$ and $\sigma_r=\tan{\delta_r}$, then the requirement that the flow in state-2 is parallel to the wedge implies that 
\begin{equation}\label{equ29}
\delta_i=-\delta_r\Longrightarrow{\frac{\sigma_i+\sigma_r}{1-{\sigma_i}{\sigma_r}}}=0.
\end{equation} 
\section{Condition for Regular Reflection}
We now consider $\beta_i$ and $\tan{\phi_i}$ as independent variables in their respective domains, and derive the condition which ensures that the regular reflection takes place in the neighborhood of the reflection point $A$.
On eliminating $M^2_1$ from (\ref{equ23}) and $(\ref{equ24})_2$, we get
\begin{equation}\label{equ30}
\beta_r=\dfrac{(\gamma+1)(1+{\beta^2_i}{\tan^2{\phi_i}})}{(\gamma+1){\beta_i}\sec^2{\phi_r}+(\gamma-1+2\tilde{b}\beta_i)({\beta^2_i}{\tan^2{\phi_i}}-\tan^2{\phi_r})}.
\end{equation} 
Substituting $\beta_r$ into $(\ref{equ26})_2$, we get
\begin{equation}\label{equ31}
\tan{\delta_r}=\dfrac{\tan{\phi_r}[2(1-\tilde{b}\beta_i)({\beta^2_i}{\tan^2{\phi_i}}-\tan^2{\phi_r})-(\gamma+1)(\beta_i-1)\sec^2{\phi_r}]}{\beta_i(\gamma+1)(1+{\beta_i}{\tan^2{\phi_i}})\sec^2{\phi_r}-2(1-\tilde{b}\beta_i)({\beta^2_i}{\tan^2{\phi_i}}-\tan^2{\phi_r})}.
\end{equation}
Using $(\ref{equ21})_2$ and (\ref{equ31}) in (\ref{equ29}), we obtain
\begin{align}\nonumber
\begin{split}
&({\beta_i}{\tan{\phi_i}}-\tan{\phi_r})[(\gamma+1)(\beta_i-1)(1+{\beta_i}{\tan^2{\phi_i}})\sec^2{\phi_r}\\&+2({\beta_i}{\tan{\phi_i}}+\tan{\phi_r})(1-\tilde{b}{\beta_i})((1+{\beta_i}{\tan^2{\phi_i}})\tan{\phi_r}-\tan{\phi_i}(\beta_i-1))]=0,
\end{split}
\end{align}
which implies that either
\begin{equation}\label{equ32}
\tan{\phi_r}={\beta_i}{\tan{\phi_i}}
\end{equation}
or
\begin{equation}\label{equ33}
\tan{\phi_r}=\dfrac{-\tan\phi_i(1+{\beta^2_i}{\tan^2{\phi_i}})(1-\tilde{b}{\beta_i})\pm\sqrt{F(\beta_i, \tan^2{\phi_i})}}{(1+{\beta_i}{\tan^2{\phi_i}})((\gamma+1-2\tilde{b})\beta_i-(\gamma-1))},
\end{equation}
with an appropriate branch, where $F$ is given by
\begin{align*}
\begin{split}
F(\beta_i, \tan^2{\phi_i})=&{\tan^2{\phi_i}}(1+{\beta^2_i}{\tan^2{\phi_i}})^2(1-\tilde{b}\beta_i)^2-(\beta_i-1)(1+\beta_i{\tan^2{\phi_i}})((\gamma+1-2\tilde{b})\beta_i\\&-(\gamma-1))((\gamma-1+2\tilde{b}{\beta_i}){\beta_i}{\tan^2{\phi_i}}+(\gamma+1)).
\end{split}
\end{align*}
Now using (\ref{equ32}) in (\ref{equ30}), we get
\begin{equation*}
\beta_r=\frac{1}{\beta_i}=\frac{\rho_0}{\rho_1}<1,
\end{equation*}
which violates entropy condition $\rho_1>\rho_0$, so it needs to be discarded. Thus the required condition follows from (\ref{equ33}) if, and only if, $F(\beta_i, \tan^2{\phi_i})$ is nonnegative, i.e.,
\begin{equation}\label{equ34}
F(\beta_i, \tan^2{\phi_i})\geq{0}
\end{equation}
Assuming this to be true, equation (\ref{equ30}) when used in (\ref{equ28}) yields,
\begin{equation}\label{equ35}
1<\dfrac{(\gamma+1)(1+{\beta^2_i}{\tan^2{\phi_i}})}{(\gamma+1){\beta_i}\sec^2{\phi_r}+(\gamma-1+2\tilde{b}\beta_i)({\beta^2_i}{\tan^2{\phi_i}}-\tan^2{\phi_r})}<\dfrac{\gamma+1}{\gamma-1+2\tilde{b}\beta_i}.
\end{equation}
The first inequality in (\ref{equ35}), implies that
\begin{equation}\label{equ36}
1+{\tan^2{\phi_r}}<\dfrac{2(1-\tilde{b}\beta_i)(1+{\beta^2_i}{\tan^2{\phi_r}})}{(\gamma+1-2\tilde{b})\beta_i-(\gamma-1)},
\end{equation}
which, in view of (\ref{equ33}), yields
\begin{equation}\label{equ36*}
(\beta_i-1)(\gamma+1)(1+{\beta_i}{\tan^2{\phi_i}})(1+{\beta^2_i}{\tan^2{\phi_i}})>0.
\end{equation} 
It may be noticed that the inequality (\ref{equ36*}) is always true; this, indeed, implies that (\ref{equ28}) follows from (\ref{equ34}). Further it follows from (\ref{equ33}) that the plus branch for $\tan{\phi_r}$ (with $\beta_i=1$) yields~ $\tan{\phi_r}=0$, whereas the minus branch yields~ $\tan{\phi_r}\Big|_{\beta_i=1}=-\tan{\phi_i}$, implying thereby that the branch with plus sign needs to be discarded since it is irrelevant here; for a valid solution, we use the minus branch for $\tan{\phi_r}$.
In order to solve the inequality (\ref{equ34}), we write $F$ in the following form
\begin{equation}\label{equ37}
F(\beta_i, \tan^2{\phi_i})= h_0+h_1{X_i}+h_2{X^2_i}+h_3{X^3_i}\equiv \widetilde{F}(\beta_i, X_i),
\end{equation} 
where
\begin{equation}\label{equ37*}
X_i=1+{\beta_i}{\tan^2{\phi_i}},
\end{equation}
$h_0=-(1-\tilde{b}\beta_i)^2{(\beta_i-1)^2}/{\beta_i}$,~~~$h_3=\beta_i(1-\tilde{b}{\beta_i})^2$,\\
$h_1={(1-\tilde{b}{\beta_i})^2(\beta_i-1)(3-1/{\beta_i})}-2(\beta_i-1)(1-\tilde{b}{\beta_i})((\gamma+1-2\tilde{b})\beta_i-(\gamma-1))$, and\\
$h_2=-((3\beta_i-2)(1-\tilde{b}{\beta_i})^2+(\beta_i-1)((\gamma+1-2\tilde{b})\beta_i-(\gamma-1))(\gamma-1+2\tilde{b}\beta_i))$.\\
In view of inequality (\ref{equ26*}) and the fact that $0\leq \tilde{b}<1$, it follows that $h_0$, $h_1$ and $h_2$ are negative, whereas $h_3$ is positive. From (\ref{equ37}), we have 
\begin{equation*}
\widetilde{F}(\beta_i, 0)= -(1-\tilde{b}\beta_i)^2{(\beta_i-1)^2}/{\beta_i}<0 ~~\text{and}~~\widetilde{F}(\beta_i, \infty)=\infty,
\end{equation*}
 implying thereby that the interval $(0, \infty)$ must contain a zero of $\tilde{F}$. Further, as the critical points of $\tilde{F}$, given by 
\begin{equation*}
{X^{\pm}_i}=\left\{{-h_2\pm({h^2_2-3{h_3}{h_1}})^{1/2}}\right\}/{3h_3},
\end{equation*}
are such that $\tilde{F}_{X_i{X_i}}>0$ at $X_i=X^{+}_i$ and $\tilde{F}_{X_i{X_i}}<0$ at $X_i=X^{-}_i$, it follows that the only positive zero of $\tilde{F}$ is given by
\begin{equation}\label{equ39}
X_i=\left\{-\frac{n}{2}+{\left(\frac{n^2}{4}+\frac{m^3}{27}\right)}^{1/2}\right\}^{1/3}+ \left\{-\frac{n}{2}-{\left(\frac{n^2}{4}+\frac{m^3}{27}\right)}^{1/2}\right\}^{1/3}\equiv{x^*},
\end{equation}
with
\begin{align}\nonumber
\begin{split}
m= &{-}\left\{(2-3\beta_i)^{2}(1-\tilde{b}\beta_i)^{4} + (\beta_i-1)^{2}(1 + \gamma(\beta_i-1) + \beta_i - 2\tilde{b}\beta_i)^{2}(\gamma - 1 + 2\tilde{b}\beta_i )^{2}\right. \\&
-\left. 2(\beta_i-1)(3\beta_i-2)(1-\tilde{b}\beta_i)^{2}(\gamma - 1 + 2\tilde{b}\beta_i)(\gamma - 1 + (2\tilde{b}-(\gamma+1))\beta_i)\right.\\&
- \left. 3(\beta_i-1)(1- \tilde{b}\beta_i)^{3}(-1 + (1+\tilde{b}+2\gamma)\beta_i + (\tilde{b} - 2(1+\gamma))\beta_i^{2})\right\}/{3\beta_i^{2}(1-\tilde{b}\beta_i)^{4}},
\end{split}
\end{align}
and
\begin{align}\nonumber
\begin{split}
n=& {-}\left\{27(\beta_i - 1)^{2}\beta_i(1-\tilde{b}\beta_i)^{6} + 9(\beta_i-1)(1-\tilde{b}\beta_i)^{3}(1-(1+\tilde{b}+2\gamma)\beta_i 
+ (2(1+\gamma) - \tilde{b})\beta_i^{2})\right. \\& 
\left.( (3\beta_i-2)(1-\tilde{b}\beta_i)^{2} + (\beta_i -1)(1+\gamma(\beta_i-1)+\beta_i - 2\tilde{b}\gamma)(\gamma-1+ 2\tilde{b}\beta_i))
+ 2(-1 + \gamma^{2}(\beta_i-1)^{2} \right. \\& 
\left. + 3\beta_i + (2\tilde{b}^{2} - 4\tilde{b} -1)\beta_i^{2} - (-2 + \tilde{b})\tilde{b}\beta_i^{3}+ 2\gamma(\beta_i-1)(1+ \tilde{b}(\beta_i-2)\beta_i))^{3}\right\}/{27\beta_i^{3}(1-\tilde{b}\beta_i)^{6}}.
\end{split}
\end{align}
Since $F(\beta_i, 0)$ = $\widetilde{F}(\beta_i, 1) < 0$ and $F(\beta_i, \infty)$ = $\infty$, the only positive zero of $F(\beta_i, \tan^2{\phi_i})$ is given by $X_i=x^{*}$.
In other words, the necessary condition for regular reflection to take place in the neighborhood of the reflection point $A$ follows from (\ref{equ37*}), namely,

\begin{equation}\label{equ40}
\tan^2{\phi_i}\geq(X_i-1)/{\beta_i}\equiv {\mathfrak{J}},
\end{equation}
leading us to conclude that there exists a critical ${\phi^{*}_i}$ of $\phi_i$, in the interval $(0, \pi/2)$, depending on $\beta_i$, $\tilde{b}$ and $\gamma$, given by
\begin{equation*}
\tan^2{\phi^{*}_i}={\mathfrak{J}},
\end{equation*}
such that there exists a unique state $(\rho_2, u_2, v_2, p_2)$ for each $\phi_i\geq\phi^{*}_i$ with $\phi_r$ given by  (\ref{equ33}) and satisfying the inequality (\ref{equ28}). In the absence of real gas effects $(\tilde{b}=0)$, we recover the result obtained by Chang and Chen \cite{chen1986}, who studied the problem of shock diffraction along a compressive corner.\\
The following table shows the effects of density ratio $\beta_i$ and $\tilde{b}$ on `${\mathfrak{J}}$'; this, indeed, shows that an increase in $\beta_i$ or $\tilde{b}$ causes the critical value `${\mathfrak{J}}$' to increase, implying thereby that for a regular reflection to take place, an increase either in the shock strength $\beta_i$ or in the van der Waals excluded volume requires the incident angle $\phi_i$ to be larger, relative to what it would have been in the absence of real gas effects.\\
\begin{table}[!h]
\begin{center}
     \begin{tabular}{ |l |l| l|l | l |l|l |l |l |l |}
      \hline
    $\bf {{\beta_i}_\downarrow}$ $\bf{\tilde{b}} \rightarrow$    &\bf	0&\bf	0.02&	\bf 0.04 &	\bf 0.06&	\bf 0.08&	\bf 0.1&\bf 	0.3&\bf	0.5	&\bf 0.7 \\ \hline
     \hline
\bf1.2&	0.2258&	0.2386&	0.2521&	0.2666&	0.2819&	0.2984&	0.5521&	1.2549	&6.0147\\ \hline
\bf1.4&	0.5193&	0.5456&	0.5741&	0.605	&0.6385&	0.6752&	1.3474	&4.4341	&\\ \hline
\bf1.6	&0.6975	&0.738	&0.7825	&0.8318&	0.8865&	0.9475&	2.301&	14.5824	& \\ \hline
\bf1.8	&0.8128	&0.8677	&0.9294&	0.999&	1.078&	1.1681&	3.6347&	&	 \\ \hline
\bf2	&0.89	&0.9598&	1.0398&	1.1319&	1.2387&	1.3633&	5.6841	&   & 	 \\ \hline
\bf2.2&	0.9431&	1.0281&	1.1274&	1.2442&	1.3827&	1.5483&	9.0801&	 &	 \\ \hline
\bf2.4	&0.98&	1.0805	&1.2003	&1.3444	&1.5191&	1.7329&	15.2028&&	 \\	 \hline
\bf2.6	&1.0057&	1.1221&	1.2637&	1.4377&	1.6535&	1.9242&	&&	 	 \\ \hline
\bf2.8	&1.0235	&1.1561	&1.3209&	1.5278&	1.7903&	2.1277&	&&	 	 \\ \hline
\bf3.0	&1.0357	&1.1849	&1.3742&	1.6171&	1.9327&	2.3482&	&&	 	 \\ \hline
\bf3.2	&1.0436	&1.2098	&.4252&	1.7077&	2.0831&	2.5901	& &&	 	 \\ \hline
\bf3.4	&1.0485&	1.2321&	1.4751&	1.8009&	2.244	&2.8577	& &&	 	 \\ \hline
\bf3.6&	1.0511&	1.2525&	1.5248&	1.8979&	2.4172&	3.1555&	 &&	 	 \\ \hline
\bf3.8	&1.0518	&1.2715&	1.5749&	1.9996&	2.6049&	3.4884&	&& 	 	 \\ \hline
\bf4.0	&1.0513	&1.2897&	1.6259&	2.1069&	2.8088&	3.8619&&&	\\ \hline
 \hline	 	 
 \end{tabular}
\end{center}
\caption{Values of ${\mathfrak{J}}$ influenced by $\beta$ and $\tilde{b}$.}
\end{table}
\section{Asymptotic analysis}
As pointed out earlier that in a weak shock regular reflection, it is experimentally observed that the reflected shock is no longer rectilinear, it joins the diffracted wave front $BC$ at a point $B$ that arises due to the influence of the compressive corner at $O$. Behind the reflected shock $AB$ the flow is a uniform supersonic flow; further downstream near the origin the flow is subsonic. Therefore, the state behind the reflected shock is not uniform and the system of governing equations, in self-similar coordinates,  becomes degenerate on the boundary $BD$ (see Figure $1(ii)$); so we look for an asymptotic solution to this problem, which is uniformly valid throughout the flow field. The boundary $BD$ is indeed a characteristic across which solution is continuous but discontinuities in its derivatives are permitted; the R-H conditions for the reflected shock, derived in section-3, provide the boundary conditions for the problem.
\subsection{R-H conditions for the incident shock}
We consider the incident shock with states $(\rho_0, 0, 0, p_0)$ and $(\rho_1, u_1, 0, p_1)$ on the two sides of it. Let $\epsilon>0$ be a dimensionless parameter measuring the shock strength, i.e.,
\begin{equation}\label{equ41}
\epsilon=(\rho_1-\rho_0)/{\rho_0}
\end{equation}
Then using (\ref{equ7**}), (\ref{equ2}) and (\ref{equ5*}) in (\ref{equ11}), the R-H conditions for the incident shock are given by the following relations
\begin{align}\label{equ8}
\begin{split}
&\frac{p_{1}}{p_0}=\frac{(\gamma+1)\rho_{1}-(\gamma-1)\rho_0-2\tilde{b}{\rho_{1}}}{(\gamma+1)\rho_0-(\gamma-1)\rho_{1}-2\tilde{b}{\rho_{1}}},~~~~~~~~~\rho_1>\rho_0\\&
u_1=\left(\frac{(p_1-p_0)(\rho_1-\rho_0)}{{\rho_0}{\rho}_{1}}\right)^{1/2},~~~\zeta=a_0\sec{\theta},~~~~v_1=0.
\end{split}
\end{align}
  In view of (\ref{equ41}) and (\ref{equ5*}), equations (\ref{equ8}) yield the following asymptotic expansions  of the state-1 variables as $\epsilon\rightarrow0$:
\begin{align}\label{equ42}
\begin{split}
&\frac{\rho_1}{\rho_0}=1+{\rho_{1}^{(1)}}\epsilon,\\&
\frac{p_1}{p_0}=1+{p_{1}^{(1)}}\epsilon+p_{1}^{(2)}\epsilon^2+O(\epsilon^3),\\&
\frac{\mbox{\scriptsize{U}}_1}{c_0}={\mbox{\scriptsize{U}}_{1}^{(1)}}\epsilon+{\mbox{\scriptsize{U}}_{1}^{(2)}}{\epsilon}^2+O(\epsilon^3),\\&
\frac{\mbox{\scriptsize{V}}_1}{c_0}={\mbox{\scriptsize{V}}_{1}^{(1)}}\epsilon+{\mbox{\scriptsize{V}}_{1}^{(2)}}{\epsilon}^2+O(\epsilon^3),\\&
\frac{a_1}{c_0}=\kappa_0+{a_{1}^{(1)}}\epsilon+{a_{1}^{(2)}}\epsilon^2+O(\epsilon^3),\\&
\frac{\zeta}{c_0}={\kappa_0}\sec{\theta}+\frac{{\kappa_0}(\gamma+1)\epsilon}{4(1-\tilde{b})}\sec{\theta}+O(\epsilon^2),~~~\frac{\mbox{\scriptsize{S}}_{1}-\mbox{\scriptsize{S}}_{0}}{c_v}=\frac{\gamma\epsilon^3}{12(1-\tilde{b})^3}(\gamma^2-1)+O(\epsilon^4),
\end{split}
\end{align}
where $\rho_{1}^{(1)}=1$, $p_{1}^{(1)}=\dfrac{\gamma}{(1-\tilde{b})}$, $p_{1}^{(2)}=\dfrac{\gamma(\gamma-1+2\tilde{b})}{2(1-\tilde{b})^2}$, ${\mbox{\scriptsize{U}}_{1}^{(1)}}=\kappa_0\cos\theta$, ${\mbox{\scriptsize{V}}_{1}^{(1)}}={-\kappa_0\sin\theta}$,  ${\mbox{\scriptsize{U}}_{1}^{(2)}}=\dfrac{(\gamma-3+4\tilde{b})\kappa_0\cos\theta}{4(1-\tilde{b})}$, ${\mbox{\scriptsize{V}}_{1}^{(2)}}=\dfrac{(3-\gamma-4\tilde{b})\kappa_0\sin\theta}{4(1-\tilde{b})}$, $a_{1}^{(1)}=\dfrac{\kappa_0(\gamma-1+2\tilde{b})}{2(1-\tilde{b})}$, $a_{1}^{(2)}=\dfrac{\kappa_0((\gamma-1)(\gamma-3+8\tilde{b})+8\tilde{b}^2)}{8(1-\tilde{b})^2}$, ${\kappa_0}={(1-\tilde{b})^{-(\gamma+1)/2}}$, $c_0=a_0/{\kappa_0}$~ and ~$\alpha<\theta<\pi$.\\
\subsection{R-H conditions for the reflected shock}  
We look for the asymptotic expansions of the state-2 variables in the following form
\begin{align}\label{equ44}
\begin{split}
&{\rho_2}/{\rho_0}=1+{\rho_{2}^{(1)}}\epsilon+O(\epsilon^2),\\&
{U_2}/{c_0}={U_{2}^{(1)}}{\epsilon}+O(\epsilon^2),\\&
{V_2}/{c_0}={V_{2}^{(1)}}{\epsilon}+O(\epsilon^2).
\end{split}
\end{align}
Substituting (\ref{equ7*}), (\ref{equ42}), and (\ref{equ44}) into (\ref{equ11}), we obtain the perturbed quantities as follows:
\begin{equation}\label{equ45}
{\rho_{2}^{(1)}}=2,~~{U_{2}^{(1)}}=2\kappa_0\cos\alpha\cos(\theta-\alpha),~~{V_{2}^{(1)}}=-2\kappa_0\cos\alpha\sin(\theta-\alpha).
\end{equation}
Thus, in view of (\ref{equ45}), equations (\ref{equ44}) become
\begin{align}\label{equ46}
\begin{split}
&\frac{\rho_2}{\rho_0}=1+{2}\epsilon+O(\epsilon^2),\\&
\frac{\mbox{\scriptsize{U}}_2}{c_0}={2\kappa_0\cos\alpha\cos(\theta-\alpha)}{\epsilon}+O(\epsilon^2),\\&
\frac{\mbox{\scriptsize{V}}_2}{c_0}={-2\kappa_0\cos\alpha\sin(\theta-\alpha)}{\epsilon}+O(\epsilon^2),\\&
\tan\psi=\tan(\pi/2-\alpha)+O(\epsilon),~~~~~~ \alpha<\theta<2\alpha. 
\end{split}
\end{align} 
\subsection{ First order solution in regions $\Omega_1$ and $\Omega_2$}
It follows from (\ref{equ6}), $(\ref{equ42})_1$ and $(\ref{equ46})_1$ that the solution, to the first order approximation, in regions  $\Omega_1$ and $\Omega_2$, is piecewise constant, i.e.,
\begin{equation}\label{equ50}
{\rho_i^{(1)}({\zeta},\theta)} = \left\{
  \begin{array}{l l }
    {\rho}^{(1)}_1=1, & \quad (\zeta, \theta)\in{\Omega_1},\\
    {\rho}^{(1)}_2=2, & \quad (\zeta, \theta)\in{\Omega_2}.
\end{array} \right.
\end{equation}
The solution (\ref{equ50}) shows that the point $B$, where the reflected wave merges with the diffracted wave smoothly, divides the diffracted wavefront $\zeta=a_0$ into two parts such that $\rho^{(1)}=2$ on $BD$ and $\rho^{(1)}=1$ on $BC$. 
\subsection{Asymptotic acoustic solutions in the diffracted region $\widetilde{\Omega}$}
In order to seek a uniformly valid asymptotic solution to the problem under consideration, we look for asymptotic expansions of the form
\begin{align}\label{equ47}
\begin{split}
&{\rho}/{\rho_0}=1+{\epsilon}{\tilde{\rho}^{(1)}}+{{\epsilon}^2}{\tilde{\rho}^{(2)}}+O(\epsilon^3),\\&
{\mbox{\scriptsize{U}}}/{c_0}={\epsilon}{\kappa_0}\tilde{{\mbox{\scriptsize{U}}}}^{(1)}+{{\epsilon}^2}{\kappa_0}\tilde{{\mbox{\scriptsize{U}}}}^{(2)}+O(\epsilon^3),\\&
{\mbox{\scriptsize{V}}}/{c_0}={\epsilon}{\kappa_0}\tilde{{\mbox{\scriptsize{V}}}}^{(1)}+{{\epsilon}^2}{\kappa_0}\tilde{{\mbox{\scriptsize{V}}}}^{(2)}+O(\epsilon^3),\\&
({{\mbox{\scriptsize{S}}}-{\mbox{\scriptsize{S}}_0}})/{c_v}={\epsilon}\tilde{{\mbox{\scriptsize{S}}}}^{(1)}+{{\epsilon}^2}\tilde{{\mbox{\scriptsize{S}}}}^{(2)}+O(\epsilon^3).
\end{split}
\end{align}
 Introducing the non-dimensional variable $\xi=\zeta/{c_0}$ and inserting the asymptotic expansions (\ref{equ47}) into $(\ref{equ5})_{1,2,3}$ and (\ref{equS}), we get the following system of equations for the first order perturbation variables  
\begin{align}\label{equ48*}
\begin{split}
&-{\xi^2}{\tilde{\rho}_{\xi}^{(1)}}+\kappa_0{\xi}{\tilde{{\mbox{\scriptsize{U}}}}_{\xi}^{(1)}}+\kappa_0(\tilde{{\mbox{\scriptsize{U}}}}^{(1)}+{\tilde{{\mbox{\scriptsize{V}}}}_{\theta}^{(1)}})=0,\\&
\kappa_0{\tilde{\rho}_{\xi}^{(1)}}-{\xi}{\tilde{\mbox{\scriptsize{U}}}_{\xi}^{(1)}}+(\kappa_0{(1-\tilde{b})}/{\gamma}){\tilde{\mbox{\scriptsize{S}}}_{\xi}^{(1)}}=0,\\&
-{\xi^2}{\tilde{\mbox{\scriptsize{V}}}_{\xi}^{(1)}}+\kappa_0{\tilde{\rho}_{\theta}^{(1)}}+(\kappa_0{(1-\tilde{b})}/{\gamma}){\tilde{\mbox{\scriptsize{S}}}_{\theta}^{(1)}}=0,\\&
{\tilde{\mbox{\scriptsize{S}}}_{\xi}^{(1)}}=0.
\end{split}
\end{align}
Eliminating $\tilde{{\mbox{\scriptsize{U}}}}^{(1)}$, $\tilde{{\mbox{\scriptsize{V}}}}^{(1)}$ and $\tilde{{\mbox{\scriptsize{S}}}}^{(1)}$, equations(\ref{equ48*}) yield the following equation for the unknown $\tilde{\rho}^{(1)}$
\begin{equation}\label{equ49}
{\xi^2}\left(\left({1-({\xi}/{\kappa_0})^2}\right){\tilde{\rho}_{\xi}^{(1)}}\right)_{\xi}+{\tilde{\rho}_{\theta\theta}^{(1)}}+{\xi}{\tilde{\rho}_{\xi}^{(1)}}=0,
\end{equation}
the solution of which satisfying the boundary condition (\ref{equ50}) at $\xi=\kappa_0$ and the boundary condition (\ref{equ7}) along the wedge surface in terms of the first order variables, i.e., ${\partial\tilde{{\rho}}^{(1)}}/{\partial{n}}=0$, takes the form 
\begin{align}\label{equ51}
\begin{split}
&{\tilde{\rho}^{(1)}}=1+{\frac{1}{\pi}}\arctan\left\{{\frac{(1-s^{2\mu})\cos\mu\pi}{-(1+s^{2\mu})\sin\mu\pi+2s^\mu\cos\mu\beta}}\right\}\\&~~~~~~~~~~~+{\frac{1}{\pi}}\arctan\left\{{\frac{-(1-s^{2\mu})\cos{\mu\pi}}{(1+s^{2\mu})\sin\mu\pi+2s^{\mu}\cos\mu\beta}}\right\},
\end{split}
\end{align}
where ${s}={(\xi/\kappa_0)}/{(1+\sqrt{1-(\xi/\kappa_0)^2})}\leq1$,~~$\mu={\frac{1}{2}}{\pi}/(\pi-\alpha)$, ~~$tan^{-1}:{\mathbb{R}}\rightarrow[0,\pi]$,~ and ~$\theta=\beta+\alpha$, with $\beta=0 $ on the wedge, indeed, in the limit of vanishing van der Waals excluded volume $(b=0)$, solution (\ref{equ51}) reduces exactly to the one obtained in \cite{keller1951}. The point $B$, where the reflected wave $AB$ merges into the diffracted wave tangentially, separates the front $\xi=\kappa_0$ into two parts $BC$ and $BD$, which are indeed the thin regions about the diffracted wave. It follows from (\ref{equ51}) that the derivative of the linearized solution normal to $BC$ is unbounded, whereas both the normal and tangential derivatives of the linearized solution are unbounded in the neighborhood of $B$. Therefore, we need an asymptotic expansion near the wavefront $\xi=\kappa_0$ in which nonlinear effects are significant; note that (\ref{equ49}) becomes degenerate at $\xi=\kappa_0$.   
The asymptotic behavior of (\ref{equ51}) near $\xi=\kappa_0$ is given by
\begin{equation}\label{equ52}
\tilde{{\rho}}^{(1)}=\rho_i^{(1)}({\zeta},\theta)+{\frac{1}{\pi}}\left({\frac{\sqrt{2}\mu\sin{2\mu\pi}}{\cos^2\mu\beta-\sin^2\mu\pi}}\right){\sqrt{1-\frac{\xi}{\kappa_0}}}+{O\left(1-\frac{\xi}{\kappa_0}\right)},
\end{equation}
which ceases to be valid at $B(\kappa_0,2\alpha)$ because $\sin\mu(\pi-2\alpha)=\sin\mu{\pi}$; we, therefore, need a different asymptotic expansion using the method of matched expansion. We will discuss the asymptotic expansion for $\theta=2\alpha$ in section~(5.6).\\
In view of (\ref{equ52}), the linear approximation of the solution near $s=1$ for $\theta\neq 2\alpha$, $(\ref{equ47})_1$ yields 
\begin{equation*}
{\rho}/{\rho_0}=1+{\epsilon}\rho_i^{(1)}({\zeta},\theta)+{\frac{\epsilon}{\pi}}\left({\frac{\sqrt{2}\mu\sin{2\mu\pi}}{\cos^2\mu\beta-\sin^2\mu\pi}}\right){\sqrt{1-\frac{\xi}{\kappa_0}}}+O\left(\frac{\epsilon^2}{\sqrt{1-{\xi}/{\kappa_0}}}\right),
\end{equation*}
which, using the polar form of $\xi$, can be written as
\begin{equation}\label{equ54}
{\rho}/{\rho_0}=1+{\epsilon}\rho_i^{(1)}+{\frac{\epsilon}{\pi}}\left({\frac{\sqrt{2}\mu\sin{2\mu\pi}}{\cos^2\mu\beta-\sin^2\mu\pi}}\right){\sqrt{1-\frac{r}{c_0\kappa_0t}}}+O\left(\frac{\epsilon^2}{\sqrt{1-{r}/{c_0\kappa_0t}}}\right).
\end{equation}
\subsection{ Nonlinear approximation}
Equations $(\ref{equ5})_{1,2,3}$ and (\ref{equS}) can be written in polar coordinates $(r,\theta)$, assuming the fact that the wedge is  symmetric about the line  $\theta=0$ i.e.,
\begin{align}\label{equ60*}
\begin{split}
&\rho_t+(\rho\mbox{\scriptsize{U}})_r+\frac{\rho\mbox{\scriptsize{U}}}{r}=0,\\&
(\rho\mbox{\scriptsize{U}})_t+(\rho\mbox{\scriptsize{U}}^2+p)_r+\frac{\rho(\mbox{\scriptsize{U}}^2-\mbox{\scriptsize{V}}^2)}{r}=0,\\&
(\rho\mbox{\scriptsize{V}})_t+(\rho\mbox{\scriptsize{U}}\mbox{\scriptsize{V}})_r+\frac{2\rho\mbox{\scriptsize{U}}\mbox{\scriptsize{V}}}{r}=0,\\&
\mbox{\scriptsize{S}}_t+\mbox{\scriptsize{U}}\mbox{\scriptsize{S}}_r=0,
\end{split}
\end{align}
where $\mbox{\scriptsize{U}}$ and $\mbox{\scriptsize{V}}$ are radial and rotational velocities defined in (\ref{equ5*}).
It may be remarked that in order to account for the nonlinear effects near the diffracted wavefront and sonic curve $\xi=\kappa_0$, where the singularity occurs, we need to construct a new expansion when $\xi$ is close to $\kappa_0$.\\
We consider a uniform state $\rho_{i}$, $\mbox{\scriptsize{U}}_{i}$, $\mbox{\scriptsize{V}}_i$, $\mbox{\scriptsize{S}}_{i}$ with $i=1,2$ as in (\ref{equ42}) or (\ref{equ46}), into which a small amplitude wave is propagating and, following the ideas of weakly nonlinear geometrical acoustics \cite{kel, keller1983, majda}, look for an asymptotic expansion for $\theta\neq2\alpha$ of the form:
\begin{align}\label{equ62}
\begin{split}
&{\rho}={\rho_i}+\delta{\hat{\rho}(r,\tau)}+{\delta^2}{\hat{\hat{\rho}}( r,\tau)}+{O(\epsilon^3)},\\&
{\mbox{\scriptsize{U}}}={\mbox{\scriptsize{U}}_i}+{\delta}\hat{\mbox{\scriptsize{U}}}(r,\tau)+{\delta^2}\hat{\hat{\mbox{\scriptsize{U}}}}(r, \tau)+{O(\epsilon^3)},\\&
{\mbox{\scriptsize{V}}}={\mbox{\scriptsize{V}}_i}+{\delta}\hat{\mbox{\scriptsize{V}}}(r,\tau)+{\delta^2}\hat{\hat{\mbox{\scriptsize{V}}}}(r, \tau)+{O(\epsilon^3)},\\&
{\mbox{\scriptsize{S}}}={\mbox{\scriptsize{S}}_i}+{\delta}\hat{\mbox{\scriptsize{S}}}+{\delta^2}\hat{\hat{\mbox{\scriptsize{S}}}}+O(\epsilon^3),
\end{split}
\end{align}
where $\tau={\delta}^{-1}\phi(r,t)$ is the `fast' variable with $\delta<<1$ as a measure of the wave amplitude.
Equations (\ref{equ60*}), in view of (\ref{equ62}), yield at $O(1)$ the following relations
\begin{align}\label{equ63}
\begin{split}
&\phi_t{\hat{\rho}}_{\tau}+{\rho_0}{\phi_r}{\hat{\mbox{\scriptsize{U}}}_{\tau}}=0,~~{c^2_0}{\kappa^2_0}{\phi_r}{\hat{\rho}_{\tau}}+{\rho_0}{\phi_t}{\hat{\mbox{\scriptsize{U}}}}_{\tau}+{{\rho_0}{c^2_0}{\kappa^2_0}(1-\tilde{b})/{\gamma}{c_v}}{\phi_r}{\hat{S}_{\tau}}=0,\\&
{\rho_0}{\phi_t}{\hat{\mbox{\scriptsize{V}}}}_{\tau}=0,~~\phi_t{\hat{\mbox{\scriptsize{S}}}}_{\tau}=0,
\end{split}
\end{align}
which, on using vector-matrix notation, can be written as
\begin{equation}\label{equ64}
{A}\widehat{W}_{\tau}=0,
\end{equation}
where $\widehat{W}=(\hat{\rho}, {\hat{\mbox{\scriptsize{U}}}}, {\hat{\mbox{\scriptsize{V}}}}, {\hat{\mbox{\scriptsize{S}}}})^{T}$ and $A=(A_{ij})$ is a $4\times4$ matrix with $A_{11}=A_{44}=\phi_t$, $A_{22}=A_{33}=\rho_0 \phi_t$, $A_{12}=\rho_0\phi_r$, $A_{21}=(c_0 \kappa_0)^2\phi_r$, $A_{24}=\rho_0 (1-\tilde{b})A_{21}/{\gamma c_v}$, and the remaining entries being zero.\\
For a nontrivial solution of (\ref{equ64}), we should have $det(A)=0$; this implies that the phase function  $\phi(r,t)$ satisfies the following eikonal equation
\begin{equation*}
\phi^2_t-{c^2_0}\kappa_0^2\phi^2_r=0,
\end{equation*}
which can be solved using the method of characteristics, showing thereby that the characteristics or rays are straight lines in the $(r,t)$ plane, and $\phi$ is constant along each ray. We shall label each ray with a parameter $\Theta$, which is constant along each ray $r=(t;\Theta)$. Thus, system (\ref{equ64}) yields 
\begin{equation}\label{equ65}
\widehat{W}=a(r, \tau; \Theta)R,
\end{equation}
where $R=({\rho_0}\phi_r, -\phi_t, 0, 0 )^T$ is the right null vector of $A$ and $a(r, \tau; \Theta)$ is an arbitrary scalar function, in which $\Theta$ occurs as a parameter.\\
Similarly, equations (\ref{equ60*}), to the order $O(\delta)$, yield
\begin{equation}\label{equ66}
A\widehat{\widehat{W}}_{\tau}+B\widehat{W}_{\tau}+C=0,
\end{equation}
where $\widehat{\widehat{W}}$ = $(\hat{\hat{\rho}}; \hat{\hat{\mbox{\scriptsize{U}}}}, \hat{\hat{\mbox{\scriptsize{V}}}}, \hat{\hat{\mbox{\scriptsize{S}}}})^{T}$, $C=(C_{i})$ and $ B=(B_{ij})$ are $4\times1$ and $4\times4$ matrices with $C_{1}=\rho_0{\hat{\mbox{\scriptsize{U}}}}_r+\rho_0{\hat{\mbox{\scriptsize{U}}}}/r$, $C_{2}=(c_0\kappa_0)^2\hat{\rho}_r$, $C_3=C_4=0$,  $B_{11}={\hat{\mbox{\scriptsize{U}}}}\phi_r$, $B_{12}={\hat{\rho}}\phi_r$, $B_{21}={\hat{\mbox{\scriptsize{U}}}}\phi_t+(c_0\kappa_0)^2(({\hat{\mbox{\scriptsize{S}}}}/{c_v})+(\hat{\rho}/{\rho_0})(\gamma-1+2\tilde{b})(1-\tilde{b})^{-1})$, $B_{22}={\hat{\rho}}\phi_t+2\rho_0{\hat{\mbox{\scriptsize{U}}}}\phi_r$, $B_{24}=(c_0\kappa_0)^2((\hat{\rho}/{c_v})+2\rho_0{\hat{\mbox{\scriptsize{S}}}}(1-\tilde{b})/{\gamma c^2_v})$, $B_{31}={\hat{\mbox{\scriptsize{V}}}}\phi_t$, $B_{32}=\rho_0{\hat{\mbox{\scriptsize{U}}}}\phi_r$, $B_{33}=\hat{\rho}\phi_t+\rho_0{\hat{\mbox{\scriptsize{U}}}}\phi_r$,  $B_{44}={\hat{\mbox{\scriptsize{U}}}}\phi_r$, and $B_{13}=B_{14}=B_{23}=B_{34}=B_{41}=B_{42}=B_{43}=0$. 
Contracting (\ref{equ66}) by the left null vector $L=(\phi_t, -\phi_r, 0, \rho_0(c_0\kappa_0)^2(1-\tilde{b})\phi_r^2/{\gamma c_v\phi_t})$ of $A$, and using (\ref{equ65}), and taking diffracted wavefront $\phi={c_0}{\kappa_0}t-r$, we get 
\begin{equation}\label{equ67}
a_r+\frac{(\gamma+1)}{2(1-\tilde{b})}a a_{\tau}+\frac{a}{2r}=0.
\end{equation}
It may be remarked that the above equation, in the absence of real gas effects, reduces to the cylindrical inviscid equation for an ideal gas reported in \cite{tesdall2}; this, indeed, implies that
\begin{equation}\label{equ69}
a=\Lambda(\Theta, w)r^{-1/2}, 
\end{equation}
along the characteristics given by  
\begin{equation}\label{equ68}
d{\tau}/{d{r}}=(\gamma+1)a/{2(1-\tilde{b})},
\end{equation}
where $\Lambda$ is an  arbitrary function describing the wave profile through its dependence on $\Theta$, which is constant along each ray, and $w$ is a fast variable that parametrizes the characteristic curves given by (\ref{equ68}). Using (\ref{equ69}) in (\ref{equ68}) and keeping in mind that $w$ is constant on the solution curves of (\ref{equ68}), we obtain
\begin{equation}\label{equ69*}
\tau=\frac{(\gamma+1)}{(1-\tilde{b})}\Lambda(\Theta, w)r^{1/2}+\chi(w),
\end{equation}
where $\chi(w)$ is an arbitrary function, which for convenience, can be replaced by $w=\psi/\delta$ with $\psi$ as a more convenient parameter. Subsequently, (\ref{equ69*}) can be expressed as
\begin{equation}\label{equ70}
\psi=\phi-{\delta}\Lambda(\Theta, \psi/\delta)(\gamma+1)r^{1/2}/{(1-\tilde{b})},
\end{equation}
which, for  $\delta\neq0$, gives  $\psi$ implicitly, a solution for which can be multivalued; this is, indeed, the nonlinearization technique   introduced by Landau \cite{landau} and Whitham \cite{whith}, which accounts for the nonlinear effects by changing the phase function in the linear solution.
Equation $(\ref{equ62})_1$, in view of equations (\ref{equ42}), (\ref{equ46}), (\ref{equ69}) and (\ref{equ65}),  becomes
\begin{equation}\label{equ71}
\rho=\rho_0+\epsilon{\rho_0}{\rho_i^{(1)}}+{\delta}{\Lambda}(\Theta, {\psi}/\delta)r^{-1/2}{\rho_0}\phi_r+O(\delta^2).
\end{equation} 
Now, for a uniformly valid solution, the linearized solution (\ref{equ54}) must match with nonlinear solution (\ref{equ71}) near the diffracted boundary $r=c_0{\kappa_0}t$, implying thereby 
\begin{equation}\label{equ71*}
{\Lambda}(\Theta, {\psi}/\delta)=-C(\beta)\sqrt{\tau},~~\epsilon\sqrt{\delta}=\delta \Rightarrow \delta=\epsilon^2, ~and~~ \Theta=\beta,
\end{equation}
where $C(\beta)=(\sqrt{2}\mu\sin{2\mu\pi})/(\pi({\sin^2\mu\pi-\cos^2\mu\beta}))$.\\
Now, in view of $(\ref{equ62})$, $(\ref{equ71})$ and $(\ref{equ71*})$,  the solution near nonlinear diffracted wavefront $\psi=0$ is given by
\begin{equation}\label{equ72}
\left(\begin{array}{c}{\rho}\\{\mbox{\scriptsize{U}}}\\{\mbox{\scriptsize{V}}}\\{\mbox{\scriptsize{S}}}\end{array}\right)=\left(\begin{array}{c}{\rho_i}\\{\mbox{\scriptsize{U}}_i}\\{\mbox{\scriptsize{V}}_i}\\{\mbox{\scriptsize{S}}_i}\end{array}\right)-{\epsilon}C(\beta)r^{-1/2}{\psi}^{1/2}\left(\begin{array}{c}{-\rho_0}\\{-c_0\kappa_0}\\{0}\\{0}\end{array}\right).
\end{equation}
Using $(\ref{equ71*})_1$ into (\ref{equ70}), we obtain
\begin{equation}\label{equ73}
\psi=\phi+\epsilon C(\beta)(\gamma+1)({\psi r})^{1/2}/(1-\tilde{b}),
\end{equation}
which on differentiating with respect to $\psi$, yields
\begin{equation}\label{equ73*}
1= \epsilon C(\beta)(\gamma+1)\psi^{-1/2} r^{1/2}/2(1-\tilde{b}).
\end{equation}
Equations (\ref{equ73}) and (\ref{equ73*}) imply that the diffracted wavefront is either a shock or a  rarefaction depending on the sign of $C(\beta)$.
\begin{equation*}
{C(\beta)} = \left\{
  \begin{array}{l l}
    {<0}, & {\beta}<{\alpha}~\hbox{(rarefaction)},\\
    {>0}, & \beta>{\alpha}~\hbox{(shock)}.\\
\end{array} \right.
\end{equation*}
\textbf{Case~1.~When diffracted wave is rarefaction ($C(\beta)<0$)}:\\
In this case, equation (\ref{equ73}) for $r<c_0{\kappa_0}t$ can be solved to yield a  positive root 
\begin{equation}\label{equ74}
\psi^{1/2}=\epsilon C(\beta)(\gamma+1) r^{1/2}/2(1-\tilde{b})+(c_0{\kappa_0}t-r+\epsilon^2 {C^2(\beta)}(\gamma+1)^2 r/4(1-\tilde{b})^2)^{1/2},
\end{equation} 
which on using in (\ref{equ72}) gives
%
\begin{align}\label{equ75}
\begin{split}
\begin{bmatrix}
{\rho}\\{\mbox{\scriptsize{U}}}\\{\mbox{\scriptsize{V}}}\\ {\mbox{\scriptsize{S}}}
\end{bmatrix}
=\begin{bmatrix}
{\rho_2}\\{\mbox{\scriptsize{U}}_2}\\{\mbox{\scriptsize{V}}_2}\\ {\mbox{\scriptsize{S}}_2}
\end{bmatrix}
& +\epsilon C(\beta){\delta^{1/2}}r^{-1/2}(\Pi+(c_0{\kappa_0}t-r+\Pi^2)^{1/2})
\begin{bmatrix}
{\rho_0}\\{c_0\kappa_0}\\{0}\\{0}
\end{bmatrix},
\end{split}
\end{align} 
where $\Pi={\epsilon C(\beta)(\gamma+1) r^{1/2}}/{2(1-\tilde{b})}$.\\
However for $r>c_0{\kappa_0}t$, in view of (\ref{equ46}), we get ~$\rho=\rho_2$,~$\mbox{\scriptsize{U}}=\mbox{\scriptsize{U}}_2$,~$\mbox{\scriptsize{V}}=\mbox{\scriptsize{V}}_2$, ~and~$\mbox{\scriptsize{S}}=\mbox{\scriptsize{S}}_2$. Indeed, when the diffracted wavefront $\psi=0$ is rarefaction, all the flow variables are continuous across it but the discontinuity occurs in their derivatives; in view of  (\ref{equ75}), the jump in the density gradient $\rho_r$, in radial direction, is given by 
\begin{equation}\label{equ76}
[\rho_r]=\frac{(1-\tilde{b})\rho_0}{(\gamma+1)r},
\end{equation}
which shows that an increase in $\tilde{b}$ causes the density gradient to decrease, implying thereby that the rarefaction wave becomes weaker and decays slowly as compared to the corresponding ideal gas ($\tilde{b}=0$) case (see Figure $3(a)$).\\
\begin{figure}[h]
\subfloat[Jump in density gradient]{
\includegraphics[width=1.9in]{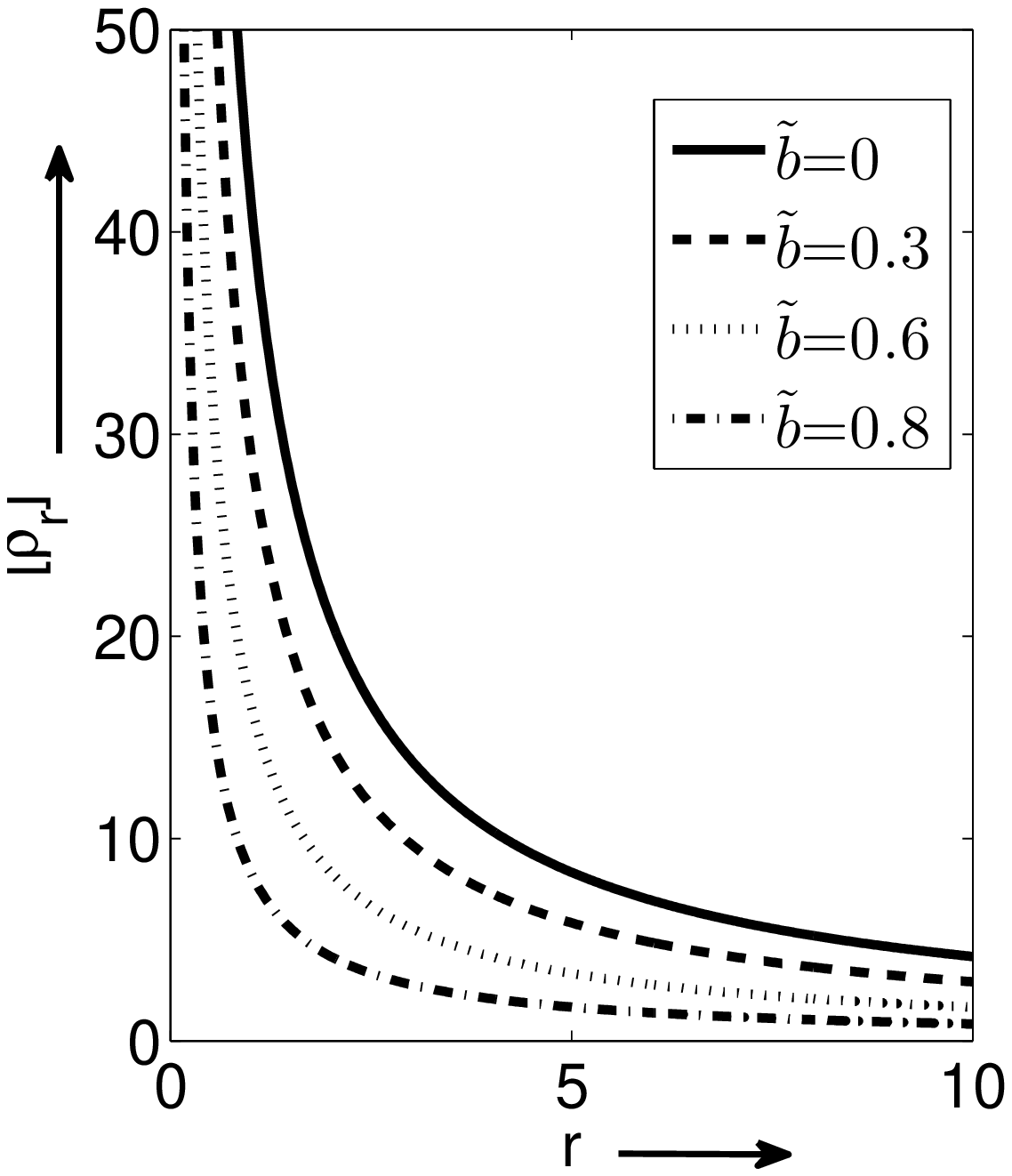}
}
\subfloat[Location]{
\includegraphics[width=2.1in]{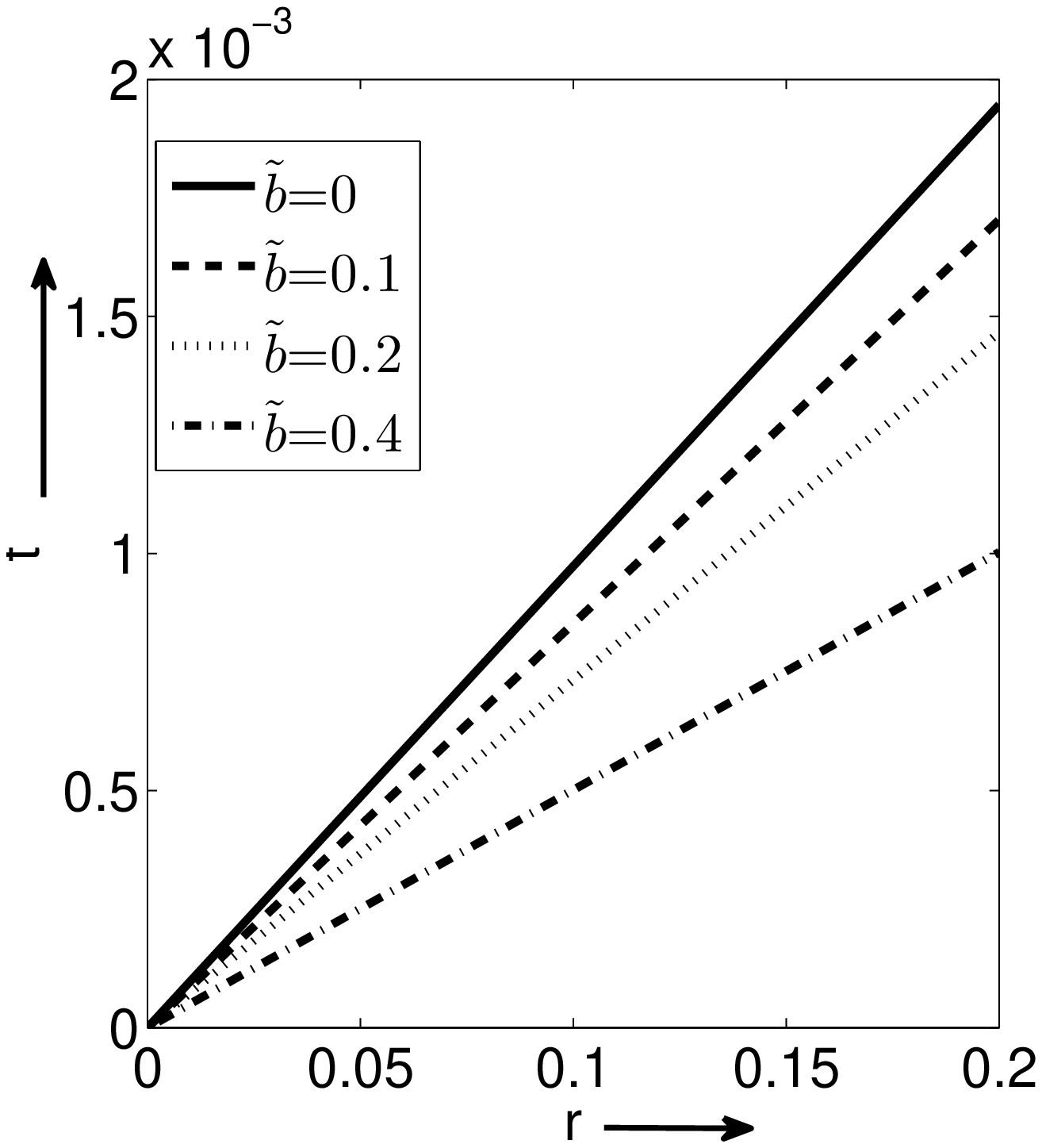}
}
\subfloat[Strength]{
\includegraphics[width=1.85in]{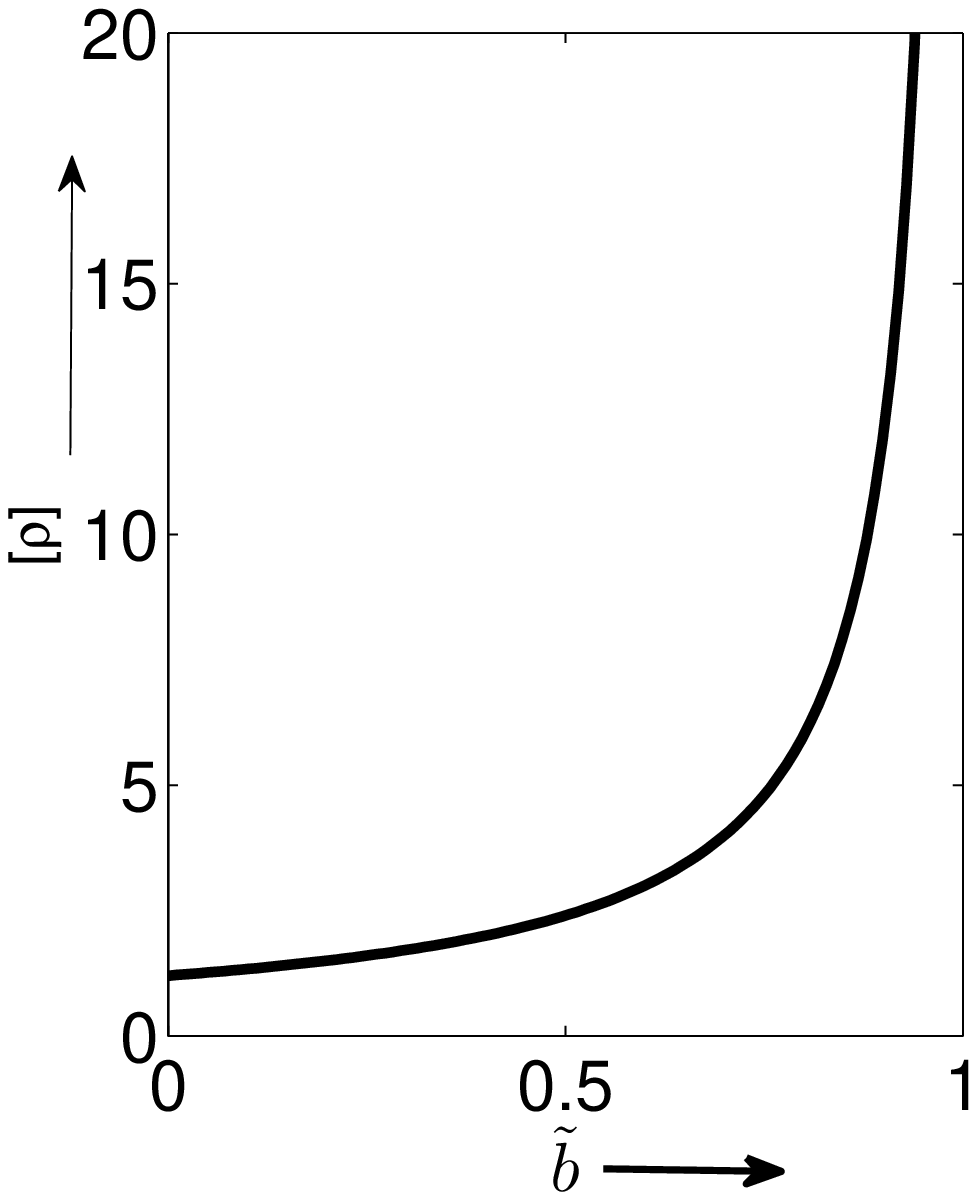}
}
\caption{\textit{(a) Jump in density gradient across the expansion wave, (b) location of diffracted shock, and (c) strength of diffracted shock for different values of $\tilde{b}$, respectively.}}
\end{figure}\\
\textbf{Case~2.~When diffracted wave is a shock ($C(\beta)>0$)}:\\ 
In this case, it is clear from (\ref{equ73*}) that the diffracted wavefront is compressive; the neighboring characteristics intersect and the envelop of the multivalued region is given by (\ref{equ73}) and (\ref{equ73*}). The multivalued region is replaced by a shock, the position of which may be determined by the equal area rule \cite{whith}. Let $r(t, \beta)$ be the shock location with $\psi_1>0$ and $\psi_2<0$ as the characteristics on the two sides of the shock. Then using the equal area rule, the equation for the shock location may be written as 
\begin{equation}\label{equ77}
r=c_0 \kappa_0 t\left(1+\frac{\epsilon^2(\gamma+1)^2C^2(\beta)}{4(1-\tilde{b})^2}\right)+O(\epsilon^4),
\end{equation}
which shows that the van der Waals excluded volume $\tilde{b}$ affects the diffracted shock location as well as the shock velocity; indeed, an increase in $\tilde{b}$ causes the velocity of diffracted shock to increase. (see Figure $3(b)$).\\
%
%
%
Now, using (\ref{equ77}) and (\ref{equ74}) in $(\ref{equ72})_1$, and matching with the boundary condition (\ref{equ41}) as $r\rightarrow c_0 \kappa_0 t$, the density immediately behind the diffracted shock is obtained as
\begin{equation}
\rho=\rho_1+ \frac{{\epsilon^2}{C^2(\beta)}(\gamma+1)}{2(1-\tilde{b})}\rho_0+O(\epsilon^3);
\end{equation}
it may be noticed that the density ahead of the shock is given by $\rho=\rho_1$. Thus, the shock strength across  the diffracted shock is given by
\begin{equation}\label{equ79}
[\rho]=\frac{{\epsilon^2}{C^2(\beta)}(\gamma+1)}{2(1-\tilde{b})}+O(\epsilon^3).
\end{equation}
From  (\ref{equ79}), it may be noticed that an increase in the van der Waals excluded volume $\tilde{b}$ causes the strength of diffracted shock to increase relative to what it would have been for the ideal gas case $(\tilde{b}=0)$ (see Figure $3(c)$).\\
\subsection{Asymptotic approximation near the singular point}
As noticed in Section (5.4) that the asymptotic behavior of the diffracted wave, given by (\ref{equ52}), breaks down in the neighborhood of $B$, we construct asymptotic expansion valid in the neighborhood of $B$ by stretching the variables $\xi$ and $\theta$; to this end we introduce new variables $r^{'}=(\xi-\kappa_0)/{\epsilon}$ and $\theta^{'}=(\theta-2\alpha)/{\epsilon^{\Delta}}$, where $\epsilon^{\Delta}$ is the gauge function with $\Delta>0$ to be determined. In terms of these new variables, the dominant part of (\ref{equ49}), after simplification, results into
\begin{equation}\label{equ52*}
2{\kappa_0}r^{'}\tilde{\rho}_{r^{'}r^{'}}^{(1)}+{\kappa_0}\tilde{\rho}_{r^{'}}^{(1)}-\epsilon^{1-2\Delta}\tilde{\rho}_{\theta^{'}\theta^{'}}^{(1)}=0.
\end{equation}
At this point, based on the principle that the leading order equation should be kept as rich as possible so that the solution contains the maximum possible information, the only choice for $\Delta$, which gives a non-degenerate reduced problem and allows all the terms in  (\ref{equ52*}) to be retained is $\Delta=1/2$. In terms of these new variables (\ref{equ51}) yields
\begin{equation}\label{equ53}
\tilde{{\rho}}^{(1)}=1+{\frac{1}{\pi}}\tan^{-1}{\frac{\sqrt{-2r'/\kappa_0}}{\theta'}}+O(\epsilon^{1/2}),
\end{equation}
where $r'<0$. \\
Accordingly, in the neighborhood of $B$, we seek asymptotic expansions of the form:
\begin{align}\label{equ80}
\begin{split}
&{\rho}={\rho_0}+{\epsilon}{\bar{\rho}(r',{\theta'})}+{\epsilon}^{3/2}{\bar{\bar{\rho}}(r',{\theta'})}+{\epsilon}^{2}{\bar{\bar{\bar{\rho}}}(r',{\theta'})}+{O(\epsilon^{5/2})},\\&
{\mbox{\scriptsize{U}}}={\epsilon}{\bar{\mbox{\scriptsize{U}}}(r',{\theta'})}+{\epsilon}^{3/2}{\bar{\bar{\mbox{\scriptsize{U}}}}(r',{\theta'})}+{\epsilon}^{2}{\bar{\bar{\bar{\mbox{\scriptsize{U}}}}}(r',{\theta'})}+{O(\epsilon^{5/2})},\\&
{\mbox{\scriptsize{V}}}={\epsilon}{\bar{\mbox{\scriptsize{V}}}(r',{\theta'})}+{\epsilon^{3/2}}{\bar{\bar{\mbox{\scriptsize{V}}}}(r',{\theta'})}+{\epsilon}^{2}{\bar{\bar{\bar{\mbox{\scriptsize{V}}}}}(r',{\theta'})}+{O(\epsilon^{3/2})},\\&
\frac{{\mbox{\scriptsize{S}}}-{\mbox{\scriptsize{S}}_0}}{c_v}={\epsilon}{\bar{\mbox{\scriptsize{S}}}(r',{\theta'})}+{\epsilon}^{3/2}{\bar{\bar{\mbox{\scriptsize{S}}}}(r',{\theta'})}+{\epsilon}^{2}{\bar{\bar{\bar{\mbox{\scriptsize{S}}}}}(r',{\theta'})}+{O(\epsilon^{5/2})}.
\end{split}
\end{align}
Substituting (\ref{equ80}) into (\ref{equ41}), and collecting respectively, $O(1)$, $O(\epsilon^{1/2})$ and $O(\epsilon)$ terms, we get
\begin{equation}\label{equ81}
{M}{\bar{W}}_{r'}=0,~~~~{M}\bar{\bar{W}}_{r'}+{N}{\bar{W}}_{\theta'}=0,~~\text{and}~~~
{M}\bar{\bar{\bar{W}}}_{r'}+{N}\bar{\bar{W}}_{\theta'}+Q=0,
\end{equation}
where $W$ is a vector of flow variables $(\rho, {\mbox{\scriptsize{U}}}, {\mbox{\scriptsize{V}}}, {\mbox{\scriptsize{S}}})^T$,\\ ${M}={\left(\begin{array}{cccc}{-c_0{\kappa_0}}&{\rho_0}&{0}&{0}\\{2{\kappa^2_0}{c_0}}&-2{\rho_0}{\kappa_0}&{0}&{\dfrac{{\rho_0}{\kappa^2_0}{c_0}}{\gamma}(1-\tilde{b})}\\{0}&{0}&-{\rho_0}{\kappa_0}&{0}\\{0}&{0}&{0}&-{\kappa_0}\end{array}\right)}$,~${N}={\left(\begin{array}{cccc}{0}&{0}&{\rho_0}/{\kappa_0}&{0}\\{0}&{0}&-{\rho_0}&{0}\\{c_0}{\kappa_0}&{0}&{0}&\dfrac{{c_0}{\kappa_0}\rho_0}{\gamma}(1-\tilde{b})\\{0}&{0}&{0}&{0}\end{array}\right)}$,\\ and~$Q=(Q_i)$ is a column vector with the components $Q_1=(\bar{\rho}\bar{\mbox{\scriptsize{U}}})_{r'}-c_0 r' \bar{\rho}_{r'}+\rho_0 \bar{\mbox{\scriptsize{U}}}/{\kappa_0}$,~ $Q_2=(\rho_0/c_0)(\bar{\mbox{\scriptsize{U}}}^2)_{r'}-2\kappa_0(\bar{\rho}\bar{\mbox{\scriptsize{U}}})_{r'}-2r'(\rho_0 \bar{\mbox{\scriptsize{U}}}-\kappa_0 c_0 \bar{\rho})_{r'}+(c_0 \kappa^2_0(\gamma-1+2\tilde{b})/{2\rho_0(1-\tilde{b})})(\bar{\rho}^2)_{r'}-\rho_0 \bar{\mbox{\scriptsize{U}}}$,~ $Q_3=(\rho_0/c_0)(\bar{\mbox{\scriptsize{U}}}\bar{\mbox{\scriptsize{V}}})_{r'}-\kappa_0(\bar{\rho}\bar{\mbox{\scriptsize{V}}}_{r'})-\rho_0 r' \bar{\mbox{\scriptsize{V}}}_{r'}$ and $Q_4=((\bar{\mbox{\scriptsize{U}}}-c_0 r')/{c_0})\bar{\mbox{\scriptsize{S}}}_{r'}$.
Let $\bar{R}$ and $\bar{L}$ be the right and left null vectors of $M$, respectively; then equation $(\ref{equ81})_1$ yields
\begin{equation}\label{equ82}
\bar{W}=U(r', \theta')\bar{R},
\end{equation}
where $\bar{R}=[\rho_0, \kappa_0 c_0, 0, 0]^T$ and $U(r', \theta')$ is an arbitrary scalar valued function.\\
In view of (\ref{equ82}), equation $(\ref{equ81})_2$ implies that
\begin{equation}\label{equ84}
\bar{\bar{W}}=V(r', \theta')(\rho_0, \kappa_0 c_0, c_0, 0)^{T},
\end{equation}
where $V(r', \theta')$ is an arbitrary scalar valued function satisfying the following relation
\begin{equation}\label{equ85}
 V_{r'}-U_{\theta'}=0.
 \end{equation}
 Now, using (\ref{equ82}) and (\ref{equ84}) in $(\ref{equ81})_3$, we obtain
 \begin{equation}\label{equ86}
 \frac{\kappa_0^2(\gamma+1)}{2(1-\tilde{b})}(U^2)_{r'}+V_{\theta'}-2\kappa_0 r'{U_{r'}}+\kappa_0 U=0.
 \end{equation}
 It may be remarked that the PDEs (\ref{equ85}) and (\ref{equ86}) bear a close structural resemblance with the self-similar (UTSD) equations analyzed in \cite{morawetz, hunter2013}.
It may be recalled that the system (\ref{equ85})-(\ref{equ86}) is the first approximation to the flow near the point $B$;
in order to obtain a uniform solution valid throughout the flow field, boundary conditions for the system (\ref{equ85})-(\ref{equ86}) must be specified in conformity with (\ref{equ42}), (\ref{equ44}), and (\ref{equ47}).
On using $\zeta=c_0(\kappa_0+\epsilon r')$ and $\theta=2\alpha+\epsilon^{1/2}\theta'$ in (\ref{equ7*}), we get the following approximations for the reflected shock in $(r', \theta')$ plane 
\begin{equation}\label{equ88}
r'=\frac{\kappa_0{\theta'}^2}{2} ~~~~~~\text{as}~~\theta'\rightarrow{\infty},
\end{equation}
showing thereby that, in the neighborhood of $B$, the straight reflected shock becomes parabolic in the limit $\theta'\rightarrow{\infty}$; we notice that an increase in the van der Waals excluded volume $\tilde{b}$ causes an increase in its latus-rectum, indicating thereby that an increase in $\tilde{b}$ causes the real gas boundaries to become larger than the corresponding ideal gas case.  \\ 
In view of (\ref{equ42}), (\ref{equ44}), (\ref{equ53}), and (\ref{equ88}), the boundary conditions for (\ref{equ85})-(\ref{equ86}) can be specified as
\begin{equation}\label{equ89}
 \lim_{{\theta'} \to \infty}{U(\eta{\kappa_0}\frac{{\theta'}^2}{2}, \theta')}\approx \lim_{{\epsilon} \to 0}{\tilde{\rho}^{(1)}}
  = \left\{
  \begin{array}{l l l}
    1, & {\eta}>1,\\
    2, & 0<{\eta}<1,\\
		{1+\frac{1}{\pi}\tan^{-1}{\sqrt{-\eta}}}, & {\eta}<0.\\
\end{array} \right.
\end{equation}
Let $r'=S_R(\theta')$ be the location of the reflected shock, which is a weak solution of the conservative system (\ref{equ85})-(\ref{equ86}). Then the jump conditions across the reflected shock may be written as 
\begin{equation}\label{equ90}
 [V]+(d{S_{R}}/{d\theta^{'}})[U]=0.
\end{equation}
\begin{equation}\label{equ91}
\frac{\kappa_0^2(\gamma+1)}{(1-\tilde{b})}[U^2/2]+(d{S_{R}}/{d\theta^{'}})[V]-2\kappa_0 S_{R}[U]=0.
\end{equation}
On using (\ref{equ90}) in (\ref{equ91}), we get
\begin{equation}\label{equ91*}
\frac{\kappa_0^2(\gamma+1)}{(1-\tilde{b})}<U>-(d{S_{R}}/{d\theta^{'}})^{2}-2\kappa_0 S_{R}=0,
\end{equation}
where $<U>$ denotes the average value of $U$ on either side of the shock.
The solution of the above equation, in view of the fact that the values of $U$ ahead and behind of $S_R$ are 1 and 2, respectively, can be written as 
\begin{equation}\label{equ92}
S_R=(\kappa_0/2)(\theta'-\theta_0)^2+(3\kappa_0/4)(\gamma+1)/(1-\tilde{b}),
\end{equation}
where $\theta_0$ is an arbitrary constant. Thus, a weak solution for the system ((\ref{equ85})-(\ref{equ86})) satisfying the boundary conditons $(\ref{equ89})_{1,2}$ can be written as
\begin{equation}\label{equ93}
 {U(r', \theta')}
  = \left\{
  \begin{array}{l l}
    1, & r'>S_R,\\
    2, & r'<S_R.\\
\end{array} \right.
\end{equation}

In a similar manner, the equation of the diffracted shock $S_D$ is obtained in the following form
\begin{equation}\label{equ94}
S_D=(\kappa_0/2)(\theta'-\theta_0)^2+(\kappa_0/4)(\gamma+1)(2+1/\pi \tan^{-1}\sqrt{-\eta})/(1-\tilde{b}),
\end{equation}
and the diffracted wave solution of the system ((\ref{equ85})-(\ref{equ86})) satisfying the boundary conditons $(\ref{equ89})_{1,3}$ can be written as
\begin{equation}\label{equ95}
 {U(r', \theta')}
  = \left\{
  \begin{array}{l l}
    1, & r'>S_D,\\
    {1+\frac{1}{\pi}\tan^{-1}{\sqrt{-\eta}}}, & r'<S_D.\\
\end{array} \right.
\end{equation}
For smooth solutions, we can eliminate $V$ from (\ref{equ85}) and (\ref{equ86}) to obtain
\begin{equation}\label{equ87}
U_{\theta'\theta'}+2\kappa_0(\vartheta U-r')U_{r'r'}+2\kappa_0 \vartheta U_{r'}^2-\kappa_0 U_{r'}=0,
\end{equation}
where $\vartheta=(\kappa_0/2)(\gamma+1)(1-\tilde{b})^{-1}$.\\
It may be noticed that equation (\ref{equ87})  is of mixed type, namely, it is hyperbolic when $\vartheta<r'$, and elliptic when $\vartheta U>r'$, however, when $\vartheta U=r'$, it corresponds to two sonic lines, $R:r'=2 \vartheta$ and $S:r'=\vartheta$. Indeed, at $r'=2\vartheta$, the reflected shock starts bending and merges asymptotically into the diffracted shock $S_D$; however the sonic line $r'=\vartheta$ is asymptotic to the diffracted shock $S_D$  in the neighborhood of the point $B'$ (see Figure $4$). \\
It may be observed that the streching transformation $r'\rightarrow {h^2}r'$, $\theta'\rightarrow {h}\theta'$, $U\rightarrow {h^2}U$, for every parameter $h>0$, leaves the equation (\ref{equ87}) invariant and, therefore, it admits a similarity solution of the form $U={\theta'}^2 f(r'/{\theta'}^2)$ such that
\begin{equation}\label{equ96}
(4x^2+(2\kappa_0/{\theta'}^2)(\vartheta{\theta'}^2f-r'))f''-(\kappa_0+2x)f'+2\kappa_0(f')^2+2f=0,
\end{equation}
with $x={r'}/{\theta'}^2$; further, as the homogeneous equation (\ref{equ96}) admits a solution of the form $f(x)=\sqrt{x}$, an expansion wave solution of (\ref{equ87}) in the region $E$ between sonic lines $r'=\vartheta$ and $r'=2\vartheta$ and satisfying the boundary conditions $(\ref{equ89})_{2,3}$ can be written as
\begin{equation}\label{equ98}
{U} = \left\{
  \begin{array}{l l l}
    {1+\frac{1}{\pi}\tan^{-1}{\sqrt{-\eta}}}, &x< \vartheta/{\theta'}^2,\\
    {{\theta'}^2 \sqrt{x}}, &\vartheta/{\theta'}^2 <x< 2\vartheta/{\theta'}^2,\\
    {2}, & x> 2\vartheta/{\theta'}^2.\\
\end{array} \right.
\end{equation}
\begin{figure}[h]
\scalebox{.65}{
\begin{tikzpicture}
			\begin{scope}
			\draw [->][thick](-4,0)--(2,0);
			 \draw [->][thick](-2.4,-3.5)--(-2.4,3.5);
	  \draw [thick](-2,0)--(-2,2.5);
	  \draw [dashed](-2,-2.9)--(-2,0);
	  \draw [thick,dashed](-2.1,-2.9)--(-2.1,2.5);
	  
\draw [thick,-][rotate around={90:(-2.4,0)}] (-5.6,3.5) parabola bend (-2.4,0) (1.2,3.5);
\draw [-][thick](-1.7,-2.9)--(-1.7,-0.9);
       \draw [rotate around={-90:(-2,0)}][thick,-] (-2,-0.5) parabola bend (-2,0)  (0.5,2);
       \draw [dashed](-1.75,-2.9)--(-1.75,-.75);
       \draw [dashed](-1.8,-2.9)--(-1.8,-.8);
       \draw [dashed](-1.85,-2.9)--(-1.85,-.7);
       \draw [dashed](-1.9,-2.9)--(-1.9,-.5);
       \draw [dashed](-1.95,-2.9)--(-1.95,-.4);
       \end{scope}
	\draw (2.1,-.2) node {$r'$};	
	\draw (-2.7,3.4) node {$\theta'$};
	\draw (-2.3,-.2) node {$B'$};
	\draw (0.5,-2.5) node {$S_R$};
	\draw (-1.6,1.6) node {$S_D$};
	\draw (-2.35,-3) node {$S$};
  \draw (-1.4,-3) node {$R$};
  \draw (-.7,0.3) node {$(1)$};
  \draw (-1.4,-1.8) node {$(2)$};
  \draw (-2,-4.5) node {$(i)$};
  \draw (-1.9,-3.3) node {$\underbrace{E}$};
   \end{tikzpicture}\hspace{-.5cm}
   \begin{tikzpicture}
			\begin{scope}
			\draw [->][thick](-4,0)--(2,0);
	  \draw [->][thick](-2.6,-3.5)--(-2.6,3.5);
   \draw [thick,-][rotate around={90:(-1.6,0)}] (-4.1,3.5) parabola bend (-1.6,0) (1.1,3.5);
   \draw [rotate around={-90:(-1.1,0)}][thick,-] (-1.1,-0.6) parabola bend (-1.1,0)  (0.7,2.3);
       \draw [thick,-](-1.1,0)--(-1.1,2.5);
       \draw [thick,dashed](-1.2,-3.2)--(-1.2,2.5);
       \draw [thick,-](0.2,-3.2)--(0.2,-1.4);
       \draw [dashed](0.1,-3.2)--(0.1,-1.35);
      \draw [dashed](0,-3.2)--(0,-1.3);
      \draw [dashed](-0.1,-3.2)--(-0.1,-1.25);
      \draw [dashed](-0.2,-3.2)--(-0.2,-1.2);
      \draw [dashed](-0.3,-3.2)--(-0.3,-1.1);
      \draw [dashed](-0.4,-3.2)--(-0.4,-1.05);
      \draw [dashed](-0.5,-3.2)--(-0.5,-1);
      \draw [dashed](-0.6,-3.2)--(-0.6,-.9);
      \draw [dashed](-0.7,-3.2)--(-0.7,-.8);
      \draw [dashed](-0.8,-3.2)--(-0.8,-.7);
      \draw [dashed](-0.9,-3.2)--(-0.9,-.6);
      \draw [dashed](-1,-3.2)--(-1,-.5);
      \draw [dashed](-1.05,-3.2)--(-1.05,-.4);
      \end{scope}
      \draw (2.1,-.2) node {$r'$};	
	\draw (-2.9,3.4) node {$\theta'$};
      \draw (-1.5,-.4) node {$B'$};
      \draw (1.5,-2.2) node {$S_R$};
      \draw (-.8,2.1) node {$S_D$};
      \draw (-1.5,-3.1) node {$S$};
  		\draw (0.5,-3) node {$R$};
  		\draw (-.5,.9) node {$(1)$};
  		\draw (.5,-2) node {$(2)$};
  		\draw (-2,-4.5) node {$(ii)$};
  		\draw (-.5,-3.5) node {$\underbrace{E}$};
  		\end{tikzpicture}\hspace{.25cm}
  		\begin{tikzpicture}
			\begin{scope}
			\draw [->][thick](-3,0)--(4,0);
	  \draw [->][thick](-2.2,-3.5)--(-2.2,3.5);
  	\draw [thick,-][rotate around={90:(0,0)}] (-2,3.5) parabola bend (0,0) (2,3.5);	
  	\draw [rotate around={-90:(1.1,0)}][thick,-] (1.1,0.5) parabola bend (1.1,0) (2.5,3.5);
      \draw [thick,-](1.1,0)--(1.1,2.5);
      \draw [thick,dashed](1,-3.2)--(1,2.5);
  	 	 \draw [thick,-](3.7,-3.2)--(3.7,-1.2);
  	 	 \draw [dashed](3.5,-3.2)--(3.5,-1.1);
  	 	 \draw [dashed](3.3,-3.2)--(3.3,-1.05);
  	 	 \draw [dashed](3.1,-3.2)--(3.1,-1);
  	 	 \draw [dashed](2.9,-3.2)--(2.9,-1);
  	 	 \draw [dashed](2.7,-3.2)--(2.7,-.95);
  	 	 \draw [dashed](2.5,-3.2)--(2.5,-.9);
  	 	 \draw [dashed](2.3,-3.2)--(2.3,-.85);
  	 	 \draw [dashed](2.1,-3.2)--(2.1,-.8);
  	 	 \draw [dashed](1.9,-3.2)--(1.9,-.75);
  	 	 \draw [dashed](1.7,-3.2)--(1.7,-.6);
  	 	 \draw [dashed](1.5,-3.2)--(1.5,-.5);
  	 	 \draw [dashed](1.4,-3.2)--(1.4,-.45);
  	 	 \draw [dashed](1.3,-3.2)--(1.3,-.4);
  	 	 \draw [dashed](1.2,-3.2)--(1.2,-.3);
  	 	 \draw [dashed](1.1,-3.2)--(1.1,0);
  	 	 	\end{scope}
	\draw (4.1,-.2) node {$r'$};	
	\draw (-2.5,3.4) node {$\theta'$};
	\draw (0.2,-.2) node {$B'$};
	\draw (5,-1.6) node {$S_R$};
	\draw (1.5,2.1) node {$S_D$};
	\draw (.8,-3.1) node {$S$};
  				\draw (4,-3) node {$R$};
  				\draw (2.4,-3.5) node {$\underbrace{E}$};
  				\draw (2,.7) node {$(1)$};
  				\draw (4.1,-1.7) node {$(2)$};
  				\draw (-1,-4.5) node {$(iii)$};
  				\end{tikzpicture}
   }
\caption{\textit{Asymptotic solution in neighborhood of $B'$, which corresponds to the point $B$ in Figure $1(ii)$; $R$ and $S$ are the sonic lines. Figure $4(i)$ corresponds to a perfect gas case $(\tilde{b}=0)$, whereas $4 (ii)$ and $4 (iii)$ account for the real gas effects with $\tilde{b}=0.3$ and $\tilde{b}=0.6$, respectively, with $\gamma=1.4$.}}
  \label{figure1}
  \end{figure}
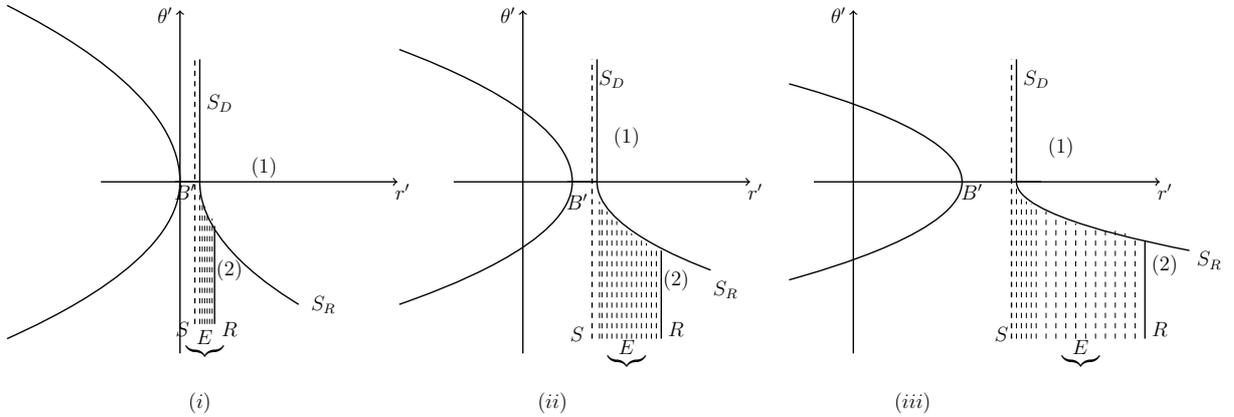 

It may be noticed that an increase in $\tilde{b}$ not only causes the sonic lines $R$ and $S$ to shift along the positive $r'$ direction, but it also increases the breadth between them (see Figure $4$). Further, equation (\ref{equ53}) shows that an increase in $\tilde{b}$ serves to reduce the jump in the derivatives of flow variables across the expansion wave $E$; also, a change in the parabolic configuration in Figure 4 leads us to reinforce our conclusion that the domain of the elliptic region in the neighborhood of the singular point $B'$ exhibits an increase with an increase in the van der Waals parameter $\tilde{b}$.
\section{Conclusions}
In this article, we explore how the real gas effects influence the self-similar solutions of the compressible Euler equations. The regular reflection configuration and the detachment criterion, influenced by the real gas effects, are studied in detail. A necessary condition is derived for the existence of regular reflection; the manner in which it is influenced by the shock strength and the van der Waals excluded volume, is clearly brought out. In the limit of vanishing van der Waals excluded volume, the ideal gas case presented in the work of Chang and Chen \cite{chen1986}, who studied the problem of shock diffraction along a compressive corner, is recovered. It is shown that for a regular reflection to take place, there exists a critical value of the angle that the shock velocity vector makes with the shock normal; it is found that an increase in the shock strength or in the van der Waals excluded volume induces an increase in the critical value, implying thereby that an increase either in the shock strength  or in the van der Waals excluded volume requires the incident angle to be larger, relative to what it would have been in the absence of real gas effects. We find that the reflected and diffracted regions as well as their boundaries, referred to as wavefronts, are significantly influenced by the real gas effects in the sense that an increase in the van der Waals excluded volume fosters an expansion of the linearized solution domain. As the state behind the reflected shock is not uniform, and the system of governing equations becomes degenerate on the boundary (referred to as the sonic arc), across which there is a continuous transition from the supersonic region to the subsonic region, we look for a uniformly valid asymptotic approximations in the flow field. Following the ideas of weakly nonlinear geometrical acoustics \cite{kel, keller1983, majda}, we construct weakly nonlinear solutions in these regions and match them with the linearized solution. It is shown that if the diffracted wave is a rarefaction wave, it gets weakened by the real gas effects and decays slowly as compared to the corresponding ideal gas case. However, when the diffracted wave is a shock, we obtain an equation for its asymptotic location, showing thereby that the real gas effects serve to enhance the speed and strength of the diffracted shock wave. Asymptotic expansions are constructed near the singular point, where the $O(\epsilon)$ approximation ceases to be valid; these expansions lead to a pair of PDEs, which bear a close structural resemblance with the self-similar (UTSD) equations, analyzed in \cite{morawetz, hunter2013}. We obtain an equation for the asymptotic position of a reflected shock in the neighborhood of the singular point, reinforcing our conclusion that the real gas effects engender the real gas boundaries to inflate. Positions of the sonic lines, at which the reflected shock starts bending and the equations change type, are determined. It is shown that the asymptotic system of coupled equations, that hold near the singular point, admits an exact similarity solution; a rarefaction wave solution, satisfying the specific boundary conditions, is obtained. It is concluded that the real gas effects serve to weaken the rarefaction wave and to enlarge the diffracted wave region supporting our earlier viewpoint.

\vspace{.2in}

\noindent
{\bf Neelam Gupta} \\
Department of Mathematics\\
Indian Institute of Technology Bombay, Powai, Mumbai-400076, India\\
E-mail: neelam@math.iitb.ac.in

\vspace{.1in}

\noindent
{\bf V. D. Sharma} \\
Department of Mathematics\\
Indian Institute of Technology Bombay, Powai, Mumbai-400076, India\\
E-mail: vsharma@math.iitb.ac.in 
\end{document}